\documentclass[twocolumn]{autart} 
\usepackage{amsfonts,amssymb,amsmath}
\usepackage{graphicx,graphics,epsfig,color}
\usepackage{rgsMacros}
\usepackage{array,amssymb,color,xspace,psfrag}
\usepackage[numbers,sort&compress,square]{natbib} 
\usepackage{amssymb,enumerate}
\usepackage{flushend}
\usepackage[dvipsnames]{xcolor}

\let\labelindent\relax
\usepackage{enumitem}
\usepackage{balance}
\usepackage{wrapfig}

\newcommand\red[1]{{\color{red}#1}}
    \definecolor{gray}{rgb}{0.33,0.4,0.47}
    \definecolor{steelblue}{rgb}{0,.42,.7}
    \definecolor{britishgreen}{rgb}{0,0.26,0.15}
    \definecolor{navyblue}{rgb}{0,0,.8}
    \definecolor{olivegreen}{rgb}{0.14,0.29,0}
    \definecolor{myred}{rgb}{0.86,0.1,0.16}

\usepackage{lipsum}

\newtheorem{definition}{Definition}
\newtheorem{lemma}{Lemma}
\newtheorem{theorem}{Theorem}
\theoremstyle{plain}

\newtheorem{proposition}{Proposition}
\newtheorem{assumption}{Assumption} 
\newtheorem{remark}{Remark}

\newtheorem{example}{Example}

\newenvironment{proof}{\noindent {\it Proof.}}{ \hfill \mbox{\footnotesize $\blacksquare$ } \\ }

\newif\ifitsdraft

\begin{document}

\begin{frontmatter}
\title{Sufficient Conditions for Robust Safety in Differential Inclusions Using Barrier Functions}

\author[gipsa]{Mohamed Maghenem},   \
\ead{mohamed.maghenem@gipsa-lab.fr}
\author[ucsc]{Masoumeh Ghanbarpour},  and \
\author[centsup]{Adnane Saoud}

\address[gipsa]{Universit\'e Grenoble Alpes, CNRS,  Grenoble-INP,  GIPSA-lab,  F-38000,  Grenoble,  France. }
\address[ucsc]{Department of Electrical and Computer Engineering, University of California, Santa Cruz, USA. }
\address[centsup]{CentraleSupelec,  University ParisSaclay,  France, and UM6P College of Computing, Benguerir, Morocco. }

\maketitle


\begin{abstract} 
This paper considers  \textit{robust-safety} notions for differential inclusions. A general framework is proposed to certify the considered 
notions in terms of barrier functions.  While most existing literature studied only what we designate by  \textit{uniform robust safety},  in this paper,  we make a clear distinction between uniform and non-uniform robust-safety.  
For both cases, we establish sufficient (infinitesimal) conditions on the nominal (unperturbed) system. That is,  our conditions involve only the barrier function and the system's right-hand side.  Our results allow for unbounded safety regions as well as non-smooth barrier functions.  Throughout the paper,  examples are provided to illustrate our results.  
\end{abstract}

\end{frontmatter}
 
\section{Introduction}    
Consider a continuous-time system of the form 
\begin{align} \label{eq.1}
\Sigma : \quad  \dot x \in F(x) \qquad x \in \mathbb{R}^n.
\end{align}
Safety is the property that requires the solutions to $\Sigma$ starting from a given set of initial conditions  $X_o \subset \mathbb{R}^n$ to never reach a given unsafe set $X_u \subset \mathbb{R}^n$
\cite{prajna2007framework,wieland2007constructive}.   Depending on the considered application, reaching the unsafe set can correspond to the impossibility of applying a predefined feedback law due the presence of saturation or,  simply,  colliding with obstacles or other systems \cite{tanner2003stable,glotfelter2017nonsmooth, 8625554}.   Guaranteeing safety is in fact a key step in many control applications (see~\cite{ames2019control} and the references therein). 

Analogous to Lyapunov 
theory for stability, barrier functions are one of the most used qualitative tools to study safety without computing the system's solutions.  Generally speaking,  a barrier function candidate is a scalar function $B : \mathbb{R}^n \rightarrow \mathbb{R}$ with a given sign on $X_o$ and the opposite sign on $X_u$.   Safety is guaranteed provided that the variation of the barrier function candidate along the right-hand side $F$,  evaluated on a particular region of the state space,  satisfies a certain inequality constraint. 
This class of sufficient (infinitesimal) conditions for safety is 
well-documented in the literature,  under different assumptions on $F$,  the  barrier function candidate $B$,  and the right-hand term of the inequality constraint; see \cite{konda2019characterizing, 10.1007/978-3-642-39799-8_17,  
FEKETA20201949,  draftautomatica,  9705088} for more details.

 When a safe system $\Sigma$ is subject to additive (possibly state-dependent) perturbation,  the classical safety property can be lost even for arbitrarily small perturbations; Example 1 in \cite{9683684}. As a consequence, the classical safety notion is not robust in nature \cite{wisniewski2016converse},  
 and a robust notion,  that can cope with additive  perturbations,  is proposed in this work. Indeed,  we  say that  $\Sigma$  is  \textit{robustly safe} with respect to the initial and unsafe sets  $(X_o,X_u)$ if we can find a continuous and strictly positive  function $\epsilon : \mathbb{R}^n \rightarrow \mathbb{R}_{>0}$,  named \textit{robustness margin},  under which, the resulting perturbed system
 \begin{align} \label{eq.2}
\Sigma_\epsilon : ~ \dot{x} \in F(x) + \epsilon(x) \mathbb{B}  \qquad  x \in \mathbb{R}^n
\end{align} 
 is safe with respect to $(X_o,X_u)$.  Here, we use $\mathbb{B} \subset \mathbb{R}^n$ to denote the closed unit ball centered at the origin.   This notion is necessarily stronger than just safety and provides a margin for safety,  which can be arbitrary small. Furthermore, when a robustly safe system $\Sigma$ admits a robustness margin $\epsilon$ that is constant, we say that $\Sigma$ is \textit{uniformly robustly safe}.   Note that uniform robust safety coincides with the robust-safety notion studied in existing literature \cite{ratschan2018converse, wisniewski2016converse,  liu2020converse, JANKOVIC2018359,seiler2021control}.   In particular, it is shown in 
\cite{wisniewski2016converse, ratschan2018converse} that,  when the complement of the unsafe set is bounded,  uniform robust safety holds if the scalar product between the gradient of $B$ and the dynamics $F$,  evaluated on the set $\partial K := \left\{ x \in \mathbb{R}^n : B(x) = 0 \right\}$,   is strictly negative. \blue{Furthermore, in \cite{liu2020converse, JANKOVIC2018359}, an upper bound on the perturbation is assumed to be known and it is incorporated in 
the barrier-function-based safety condition. Note that unbounded barrier functions are used in \cite{JANKOVIC2018359}. Finally, in \cite{seiler2021control}, 
 robust safety with respect to specific classes of unknown dynamics, whose solutions verify a certain integral constraint is considered. }

In this paper,  we consider a general notion of robust safety in the context continuous-time systems modeled by differential inclusions.  In particular,  we put the proposed notion in perspective with respect to existing literature, as we make a clear distinction between uniform and non-uniform robust safety. 
\blue{After that, we motivate the utility of the considered notions on three different contexts reminiscent of self-triggered control, fast-slow dynamics, and assume-guarantee contract}. Furthermore, we propose sufficient (infinitesimal) conditions, involving only the barrier function candidate $B$ and the nominal right-hand side $F$ to certify the two robust-safety notions. In comparison to the existing works in  \cite{wisniewski2016converse,ratschan2018converse,liu2020converse}, our contributions are as \blue{follows}.  First,  we consider the general context of differential inclusions while only differential equations are considered in existing works.  Second,  we allow for state-dependent robustness margins and our results do not require a bounded safety region $\mathbb{R}^n \backslash X_u$.    Finally,  we allow for \blue{barrier functions that are not necessarily smooth}.  In particular, to certify robust safety,  we allow for barrier functions that are either continuously differentiable or only locally Lipschitz.  Moreover,  to certify uniform robust safety,  we allow for barrier functions that are either continuously differentiable,  locally Lipschitz, or only semicontinuous. 

\blue{The remainder of this paper is organized as follows. Preliminaries on nonsmooth analysis, set-valued maps,  and differential inclusions are in Section \ref{sec.2}. Safety and the considered robust-safety notions are in Section \ref{sec.not}. Motivations of the proposed robust-safety notions are in Section \ref{sec.3-}. The main results are in Section \ref{sec.3}. Some numerical examples are in Section \ref{sec.exp}. The paper is concluded with some perspectives and future works.

Different from the conference version \cite{9683684},  here we include proofs, explanations, and  motivational examples.} 

\textbf{Notations.}  For $x \in \mathbb{R}^n$ and $y \in \mathbb{R}^n$, $x^{\top}$ denotes the transpose of $x$, $|x|$ the Euclidean norm of $x$ and $\langle x, y \rangle:= x^\top y $ the inner product between $x$ and $y$.  For a set $K \subset \mathbb{R}^n$, we use $\mbox{int}(K)$ to denote its interior, $\partial K$ to denote its boundary,   $U(K)$ an open neighborhood around $K$,  \blue{and $|x|_K$ to denote the distance between $x$ and the set $K$}.     For $O \subset \mathbb{R}^n$,  $K \backslash O$ denotes the subset of elements of $K$ that are not in $O$.    
For a function $\phi : \mathbb{R}^n \rightarrow \mathbb{R}^m$, $\dom \phi$ denotes the domain of definition of $\phi$.  By $F : \mathbb{R}^m \rightrightarrows \mathbb{R}^n $, we denote a set-valued map associating each element $x \in \mathbb{R}^m$ into a subset $F(x) \subset \mathbb{R}^n$.  For a set $D \subset \mathbb{R}^m$,  $F(D) := \{ \eta \in F(x) : x \in D \}$.    For a differentiable map $x \mapsto B(x) \in \mathbb{R}$,  $\nabla_{x_i} B$ denotes the gradient of $B$ with respect to $x_i$,  $i \in \{1,2,...,n\}$, and \blue{$\nabla B$ denotes the gradient of $B$. 
}

\section{Preliminaries} \label{sec.2}
\subsection{Set-Valued vs Single-Valued Maps} \label{sec.1a}
Given a \blue{scalar} map $B: K \rightarrow \mathbb{R}$, where $K \subset \mathbb{R}^m$. 
\begin{itemize}
\item $B$ is said to be \textit{locally Lipschitz} if,  for each nonempty compact set  $U \subset K$,  there exists $k>0$ such that, 
for each $(x_1,x_2) \in (U \cap K) \times (U \cap K)$, 
$|B(x_1) - B(x_2)| \leq k |x_1-x_2|$.  
\item  $B$ is said to be \textit{lower semicontinuous} at $x \in K$ if, for every sequence $\left\{ x_i \right\}_{i=0}^{\infty} \subset K$ such that $\lim_{i \rightarrow \infty} x_i = x$, we have $\liminf_{i \rightarrow \infty} B(x_i) \geq B(x)$. 
\item  $B$ is said to be \textit{upper semicontinuous} at $x \in K$ if, for every sequence $\left\{ x_i \right\}_{i=0}^{\infty} \subset K$ such that $\lim_{i \rightarrow \infty} x_i = x$, we have $\limsup_{i \rightarrow \infty} B(x_i) \leq B(x)$.  
\item  $B$ is said to be \textit{continuous} at $x \in K$ if it is both upper and lower semicontinuous at $x$. 
\end{itemize}

$B$ is said to be upper, lower semicontinuous, or continuous if, respectively, it is upper, lower semicontinuous, or continuous for all $x \in K$.

Given a set-valued map $F: K \rightrightarrows \mathbb{R}^n$, where $K \subset \mathbb{R}^m$. 

\begin{itemize}
\item $F$ is said to be \textit{locally Lipschitz} if,  for each nonempty compact set $U \subset K$, there exists $k>0$ such that, 
for all $(x_1,x_2) \in (U \cap K) \times (U \cap K)$, 
$F(x_1) \subset   F(x_2) + k |x_1-x_2| \mathbb{B}$.  
\item The map $F$ is said to be \textit{lower semicontinuous}  (or, equivalently, \textit{inner semicontinuous}) at $x \in K$ if for each 
$\varepsilon > 0$ and $y_x \in F(x)$, there exists $U(x)$ satisfying the following property: for each $z \in U(x) \cap K$, there exists $y_z \in F(z)$ such that $|y_z - y_x| \leq \varepsilon$; see  \cite[Proposition 2.1]{michael1956continuous}.  
\item The map $F$ is said to be \textit{upper semicontinuous} at $x \in K$ if, for each 
$\varepsilon > 0$, there exists $U(x)$ such that for each $y \in U(x) \cap K$, $F(y) \subset F(x) + \varepsilon \mathbb{B}$; see \cite[Definition 1.4.1]{aubin2009set}.
\item The map $F$ is said to be \textit{continuous} at $x \in K$ if it is both upper and lower semicontinuous at $x$.

\item The map $F$ is said to be \textit{locally bounded} at $x \in K$,  if there exist $U(x)$ and $\beta > 0$ such that  $|\zeta| \leq \beta$ for all $\zeta \in F(y)$ and for all $y \in U(x) \cap K$.  
\end{itemize}
The map $F$ is said to be upper, lower 
semicontinuous,  continuous, or locally bounded if, respectively, it is upper, lower 
semicontinuous,  continuous, or locally bounded for all $x \in K$.

\subsection{Tools for Nonsmooth Analysis} \label{sec.1b}

In this section, we recall some tools from \cite{clarke2008nonsmooth} that will allow us to use  non-smooth barrier functions.

\begin{definition} [Clarke generalized gradient] \label{defgen}
\blue{Let $B : \mathbb{R}^n \rightarrow \mathbb{R}$ be locally Lipschitz.  Let $\Omega \subset \mathbb{R}^n$ be any null-measure set,  and let $\Omega_B \subset \mathbb{R}^n$ be the set of points at which $B$ fails to be differentiable. The Clarke generalized gradient of $B$ at $x$ is the set-valued map 
$\partial_C B : \mathbb{R}^n \rightrightarrows \mathbb{R}^n$  given by
\begin{align*}
\partial_C B(x) := \co \left\{ \lim_{i \rightarrow \infty} \nabla B(x_i) : x_i \rightarrow x,~x_i \in \Omega \backslash \Omega_B \right\}.
\end{align*}
where $\co(\cdot)$ is the convex hull of the elements in $(\cdot)$}.  
\end{definition} 

\begin{remark}
When $B$ is locally Lipschitz, then $\partial_C B$ is upper semicontinuous and its images are nonempty, compact, and convex. 
\end{remark}

\begin{definition}[Proximal subdifferential]  \label{defps}
The proximal subdifferential of $B : \mathbb{R}^n \rightarrow \mathbb{R}$ at $x$ is the set-valued map 
$\partial_P B : \mathbb{R}^n \rightrightarrows \mathbb{R}^n$ given by
\begin{align*}
\partial_P B(x) := \left\{ \zeta \in \mathbb{R}^n : [\zeta^\top~-1]^\top \in N^P_{\epi B} (x, B(x)) \right\},
\end{align*}
where $ \epi B := \left\{(x, r) \in \mathbb{R}^n \times \mathbb{R} : r \geq B(x) \right\}$ \blue{and,  given  a closed subset $S \subset \mathbb{R}^{n+1}$ and $y \in \mathbb{R}^{n+1}$, 
$$ N_S^P(y) := \left\{ \zeta \in \mathbb{R}^{n+1} : \exists r > 0~\mbox{s.t.}~|y+r\zeta|_{S} = r 
|\zeta| \right\}. $$}
\end{definition}
\begin{remark}
\blue{When $B$ is twice continuously differentiable at 
$x \in \mathbb{R}^n$,    then 
$\partial_P B(x) = \left\{ \nabla B(x) \right\}$.
Moreover,  the latter equality is also true  if $B$ is continuously differentiable at $x$ and $\partial_P B(x) \neq \emptyset$. }
\end{remark}

\subsection{Differential Inclusions}  \label{sec.03}

\begin{definition} [Concept of solutions] \label{def.ConSol}
{Consider the map} $\phi : \dom \phi \rightarrow \mathbb{R}^n$,  
 \blue{where $\dom \phi$ is either of the form $\dom \phi = [0,T]$ for $T \in \mathbb{R}_{\geq 0}$  or of the form $\dom \phi = [0,T)$ for $T \in \mathbb{R}_{\geq 0} \cup \{+ \infty\}$}. {$\phi$}  is a solution to $\Sigma$ starting from $x_o \in \mathbb{R}^n$ if $\phi(0) = x_o$,  the map $t \mapsto \phi(t)$ is locally absolutely continuous,  and
$\dot{\phi}(t) \in F(\phi(t))$ for almost all $t \in \dom \phi$.
\end{definition}

 A solution $\phi$ to $\Sigma$ starting from $x_o \in \mathbb{R}^n$ is maximal if there is no solution $\psi$ to $\Sigma$ starting from $x_o$ such that $\psi(t) = \phi(t)$ for all $t \in \dom \phi$ and $\dom \phi$ strictly included in $\dom \psi$. 

\begin{remark}
\blue{We say that $\phi : \dom \phi \rightarrow \mathbb{R}^n$ is a backward solution to $\Sigma$ if 
 either $\dom \phi = [-T,0]$ for $T \in \mathbb{R}_{\geq 0}$  or $\dom \phi = (-T,0]$ for $T \in \mathbb{R}_{\geq 0} \cup \{+ \infty\}$  and $\phi(-(\cdot)) : - \dom \phi \rightarrow \mathbb{R}^n$ is a solution to $\Sigma^- : \dot{x} \in -F(x)$. }
 \end{remark}

The differential inclusion $\Sigma$ is said to be well posed if the right-hand side $F$ satisfies the following assumption.

\begin{assumption} \label{ass1} 
The map $F: \mathbb{R}^n \rightrightarrows \mathbb{R}^n $ is upper semicontinuous and $F(x)$ is nonempty, compact, and convex for all $x \in \mathbb{R}^n$.
\end{assumption}

\begin{remark}
Assumption \ref{ass1} guarantees the existence of solutions and adequate structural properties for the set of solutions to $\Sigma$;  see \cite{aubin2012differential, Aubin:1991:VT:120830,4518905}.   Furthermore,  when $F$ is single valued, Assumption \ref{ass1} reduces to $F$ being continuous.
\end{remark}

\section{Safety vs Robust-Safety Notions} \label{sec.not}

Given a set of initial conditions $X_o \subset \mathbb{R}^n$ and an unsafe set $X_u \subset \mathbb{R}^n$,  such that $X_o \cap X_u = \emptyset$, we recall that $\Sigma$ is safe with respect to $(X_o,X_u)$ if, for each solution $\phi$ starting from $X_o$,  we have 
$\phi(\dom \phi) \subset \mathbb{R}^n \backslash X_u$  \cite{prajna2005optimization}.

Note that safety with respect to  $(X_o, X_u)$ is verified if and only if there exists $K \subset \mathbb{R}^n$, with $X_o \subset K$ and $K \cap X_u = \emptyset$,  that is \textit{forward invariant};  namely,  for each solution $\phi$ to $\Sigma$ starting from $K$, $\phi(\dom \phi) \subset  K$. 

The latter set $K$ is usually chosen or designed to be the zero-sublevel set of a scalar function named \textit{barrier function candidate}. 
\begin{definition}[Barrier function candidate]
A scalar function $B : \mathbb{R}^n \rightarrow \mathbb{R}$ is a barrier function candidate with respect to $(X_o,X_u)$ if
\begin{equation}
\label{eq.3}
\begin{aligned} 
 B(x) > 0 \quad \forall x \in X_u \quad \text{and} \quad 
 B(x) \leq 0 \quad \forall x \in  X_o.  
\end{aligned}
\end{equation}
\end{definition}
A barrier candidate $B$ defines the zero-sublevel set 
\begin{align} \label{eq.4} 
K := \left\{ x \in \mathbb{R}^n : B(x) \leq 0 \right\},
\end{align} 
which is closed provided that $B$ is at least lower semicontinuous.  Note also  that $X_o \subset K$ and $K \cap X_u = \emptyset$.   Hence,  forward invariance of the set $K$ in \eqref{eq.4} implies safety with respect to $(X_o,X_u)$, and $B$ becomes a  \textit{barrier certificate}.

We now introduce the considered robust-safety notions.
\begin{definition} [Robust safety] \label{def-RobSafe} 
System $\Sigma$ in \eqref{eq.1} is robustly safe with respect to 
$(X_o, X_u)$ if there exists a continuous function $\epsilon : \mathbb{R}^n \rightarrow \mathbb{R}_{>0}$ such that $\Sigma_{\epsilon}$ in \eqref{eq.2} is safe with respect to $(X_o,X_u)$. Such a function $\epsilon$ is called 
\textit{robust-safety margin}.  
\end{definition} 
When $\epsilon$ is a robust-safety margin,  then every continuous function $\epsilon_1 : \mathbb{R}^n \rightarrow \mathbb{R}_{>0}$ satisfying  $\epsilon_1(x) \leq \epsilon(x)$  for all $x \in \mathbb{R}^n$ 
 is also a robust-safety margin.   

Clearly, robust safety implies  safety.   However,  the opposite is not always true,  see Example 1 in \cite{9683684}.  
 
Next,  we relate the proposed robut-safety notion to the one studied in  \cite{wisniewski2016converse,  ratschan2018converse,  liu2020converse}.  That is,  we introduce a stronger notions that requires a constant robust-safety margin.  

\begin{definition} [Uniform robust safety] \label{def-UnifRobSafe} 
System $\Sigma$ in \eqref{eq.1} is uniformly robustly safe with respect to $(X_o, X_u)$ if it is robustly safe with respect to  $(X_o, X_u)$ and admits a constant robust-safety margin.
\end{definition}

Uniform robust safety implies robust safety. The opposite is not always true,  see Example 2 in \cite{9683684}.  

\section{Motivation} \label{sec.3-}

We motivate the  considered notions via three examples.   
\begin{example}[Safety via self-triggered control] \label{sec.3-a}
\blue{
Consider the control system 
\begin{align} \label{dif.incu}
\Sigma_u :  \dot{x} = f(x,u) \qquad  x \in  \mathbb{R}^n, \quad u \in \mathbb{R}^{m}.
\end{align}   
Consider a feedback law $\kappa: \mathbb{R}^n \mapsto \mathbb{R}^{m}$  rendering the closed-loop system
$\Sigma : \dot{x} =  F(x) := f(x,\kappa(x))$ safe with respect to $(X_o,X_u) \subset \mathbb{R}^n \times \mathbb{R}^n$.  

In a self-triggered setting \cite{di_benedetto_digital_2013},  measurements of the state  are available only at  a sequence of times  
$\{t_i\}^{\infty}_{i=0} \subset \mathbb{R}_{\geq 0}$ ( $t_0 = 0$ and $t_{i+1} > t_i$) and the control law $\kappa$ remains constant between each two samples
$t_i$ and $t_{i+1}$, that is, 
\begin{align*} 
u(t) = \kappa(x(t_i)) \quad \forall t \in [t_i, t_{i+1}) \quad \forall i \in \mathbb{N}.
\end{align*}
As a result,  a solution $\phi$ to the self-triggered closed-loop system is governed by 
\begin{align} \label{eq:closedloopSys+}
   \dot{\phi}(t) =  F(\phi(t))  + \gamma(\phi(t),\phi(t_i))   \quad \forall t \in [t_i, t_{i+1}),
\end{align}
where 
$ \gamma(x,y) :=  f(x,\kappa(y)) - f(x,\kappa(x)).  $ 

The objective, in this case, is to design an algorithm that uses $\phi(t_i)$ to deduce the largest possible $t_{i+1}$ such that the self-triggered closed-loop system is safe.  Moreover,  it is important to guarantee that
$ t_{i+1} - t_i  > T > 0$ for all $i \in \mathbb{N}$, 
to exclude Zeno behaviors. 

Assuming $f$ and $\kappa$ to be  sufficiently smooth, the problem can be solved by a finding $\epsilon: \mathbb{R}^n \rightarrow \mathbb{R}_{> 0}$ continuous and a sequence $\{t_i\}^\infty_{i=1}$ such that
\begin{itemize}
\item [St1.] For each $t \in [t_i,t_{i+1})$, 
$\gamma(\phi(t),\phi(t_i))   \subset  \epsilon(\phi(t)) \mathbb{B}$.  
\item[St2.] The function $\epsilon$ is a robust-safety margin for $\Sigma$. 
\end{itemize}  
We show in Example \ref{exp...} that safety for 
$\Sigma$ is not enough to guarantee safety for
the self-triggered closed-loop system. However, robust safety can be enough, under some extra assumptions. }
\end{example}

\begin{example}[Safety under singular perturbations] \label{expsingpert}
\blue{
Consider the 
singularly-perturbed system 
\begin{equation} \label{eqsingper}
\begin{aligned} 
\dot{x} & =  f(x,e) \\
\varepsilon \dot{e} & =  g(x,e,\varepsilon)   \qquad (\varepsilon, x,e) \in \mathbb{R}_{> 0} \times  \mathbb{R}^{n} \times \mathbb{R}^{m}.
\end{aligned}
\end{equation}
  Assume that 
\begin{itemize}    
\item $e = 0$ is solution to $0=g(x,e,0)$  for all $x \notin X_u$. 
\item The system $ \dot{x} = F(x) := f( x , 0)$ is safe with respect to some $(X_o,X_u) \subset \mathbb{R}^{n} \times \mathbb{R}^{n}$. 
\end{itemize}
The objective,  in this case,   is to show that \eqref{eqsingper} is safe with respect to $(X_o \times \mathbb{B},  X_u \times \mathbb{R}^m)$  for sufficiently small values of $\varepsilon > 0$.  To do so,  we re-express \eqref{eqsingper} as  
\begin{equation} \label{eqsingper1}
\begin{aligned} 
\dot{x} & =  F(x) + \gamma_1(x,e)   \\
\varepsilon \dot{e} & = g(x,e,0) +  \gamma_2 (x,e,\varepsilon), 
\end{aligned}
\end{equation}
where $\gamma_1 (x,e) :=  f(x, e) - f(x,0)$
and $\gamma_2 (x,e,\varepsilon)  :=  g(x,e,\varepsilon) - g(x,e,0)$. 

Under some mild regularity assumptions on $f$ and $g$,  we can find $\Gamma_{1,2} :  \mathbb{R}^n \times \mathbb{R}^m \rightarrow \mathbb{R}_{\geq 0}$ continuous such that
$ \gamma_1(x,e)  \subset \Gamma_1(x,e) \mathbb{B}$,   $\Gamma_1(x ,  0) = 0$,  and
$\gamma_2(x,e,\varepsilon)  \subset \varepsilon \Gamma_2(x , e) \mathbb{B}$.  Assume further the existence of
$\Gamma_{3} :   \mathbb{R}^m \rightarrow \mathbb{R}_{\geq 0}$  continuous such that, for each $\varepsilon > 0$, we have 
\begin{align} \label{eqgamma3} 
\gamma_2(x , e, \varepsilon) \subset \varepsilon \Gamma_3(e) \mathbb{B} \qquad \forall x \notin X_u.
\end{align}
As a result,  to show safety 
for \eqref{eqsingper} with respect to 
$(X_o \times \mathbb{B}, X_u \times \mathbb{R}^m)$,  it is enough to show that the same property holds for  the differential inclusion 
\begin{equation}
\label{eqsingper3}
\begin{aligned} 
\dot{x} & \in  F(x) + \Gamma_1(x,e) \mathbb{B}  
  \\
\varepsilon \dot{e} & \in g( \mathbb{R}^n \backslash X_u, e,0) +  \varepsilon \Gamma_3(e) \mathbb{B}.  
\end{aligned}
\end{equation} 
To do so, we propose to find $R$,  $r$, and $T>0$ such that
\begin{itemize}
\item[St1.] The system $$\Sigma_{\epsilon_r} : \dot{x} \in F(x) + \epsilon_r(x) \mathbb{B},  \quad \epsilon_r(x) :=  \sup \{ \Gamma_1(x,r \mathbb{B}) \}, $$
is safe with respect $(Reach(T, X_o), X_u)$,  where $Reach(T, X_o)$ is the reachable set from $X_o$ over $[0,T]$ along the solutions to
$ \dot{x}  \in  F(x) + \Gamma_1(x,R \mathbb{B}) \mathbb{B}.    
 $

\item[St2.]   $\exists \varepsilon^*>0$  such that,  for each $\varepsilon \in (0, \varepsilon^*]$,  the solutions to the second equation in \eqref{eqsingper3},  starting from $\mathbb{B}$,  satisfies 
$ |e(t)|_{\infty} \leq R$ and  $\sup_{t \geq T} \{  |e(t)|  \}  \leq r$. 
\end{itemize}

To guarantee St1., we can use robust safety for slow system $\dot{x} = F(x)$ and make $r$ sufficiently small so that $\epsilon_r$ is sufficiently small; see Example \ref{expsing} in Section \ref{sec.exp}.  }
\end{example}

\begin{example}[Contracts via robust safety]
\label{expagcontracts}
\blue{Consider the perturbed system 
\begin{align} 
\label{eqn:sys_pert}
\Sigma_w : \quad  \dot x = f(x,w) \qquad x \in \mathbb{R}^n, \qquad w \in \mathbb{R}^m.
\end{align}
 We recall from  \cite{saoud2018composition,saoud2021assume} that an assume-guarantee contract for $\Sigma_w$
is a tuple $\mathcal C := (A,K)$, where  $A \subseteq \mathbb{R}^m$ is a set of assumptions and $K \subseteq \mathbb{R}^n$ is a closed set of guarantees.
We say that $\Sigma_w$ weakly satisfies $\mathcal C$, if, for any maximal solution pair $(w,\phi) : \dom \phi \rightarrow \mathbb{R}^m\times \mathbb{R}^n$ to $\Sigma_w$, we have
\begin{itemize}
\item $\phi(0)\in K$;
\item  $\forall t \in \dom \phi$ such that $w([0,t])\subset A$, 
$\phi([0,t]) \subset K$.
\end{itemize}
We say that $\Sigma_w$ strongly satisfies $\mathcal C$, if for any maximal solution pair $(w,\phi) : \dom \phi \rightarrow \mathbb{R}^m\times \mathbb{R}^n$ to $\Sigma_w$, we have	
\begin{itemize}
\item $\phi(0)\in K$;
\item $\forall t \in \dom \phi$ such that $w([0,t])\subset A$, there exists $\delta>0$ such that $\phi([0,t+\delta] \cap \dom \phi) \subset K$.
\end{itemize}
An assume-guarantee contract for $\Sigma_w$ states that, when the disturbance $w$ belongs to $A$ up to  $t \in \mathbb{R}_{\geq 0}$, then the state $x$ belongs to $K$ at least until $t$, or until $t+\delta$, with $\delta>0$, in the case of strong satisfaction.

Given a system composed of interconnected components, the notion of assume-guarantee contracts allows the verification of desired properties \textit{compositionally}. 
In particular, to guarantee that the system is safe, we decompose this property into local 
assume-guarantee contracts, and then verify that each component satisfies its contract. 
It is shown in \cite{saoud2018composition} that while the concept of weak contract satisfaction can be applied to cascaded compositions, the concept of strong satisfaction is crucial in the context of feedback composition; see Example \ref{examp:agc} in Section \ref{sec.exp}. 

As we show in Proposition \ref{propagc} in the Appendix, when 
\begin{align} \label{eqSigmaagc}
\Sigma : \dot x \in F(x) := \co \{f(x,A)\} \qquad x \in \mathbb{R}^n 
\end{align}
is safe with respect to $(K,\mathbb{R}^n \backslash K)$, then  $\Sigma_w$ weakly satisfies $\mathcal{C} := (A,K)$. Moreover, if $K$ is closed, $f$ is locally Lipschitz 
, and $\Sigma$ is robustly safe with respect to $(K,\mathbb{R}^n \backslash K)$, then  $\Sigma_w$ strongly satisfies  $\mathcal{C} := (A,K)$.}
\end{example}

\section{Main Results}  \label{sec.3}
\blue{Since robust safety for $\Sigma$ coincides with safety for  $\Sigma_\epsilon$,  when $\epsilon$ is a valid robustness margin, sufficient conditions for safety,  applied to  
$\Sigma_\epsilon$, can be used to certify robust safety.    
Inspired by \cite{liu2020converse},  we  formulate the following proposition, whose proof is a straightforward consequence of Lemmas \ref{lemA8} and \ref{lemA9} in the Appendix. }

\begin{proposition}
\blue{Consider system $\Sigma$ such that Assumption \ref{ass1} holds. Let $B : \mathbb{R}^n \rightarrow \mathbb{R}$ be a continuous barrier function candidate with respect to $(X_o,X_u) \subset \mathbb{R}^n \times \mathbb{R}^n$.  Then,  $\Sigma$ is robustly safe with respect to $(X_o,X_u)$ if there exists a continuous function $\epsilon : \mathbb{R}^n \rightarrow \mathbb{R}_{>0}$ such that 

\begin{enumerate}[label={($\star''$)},leftmargin=*]
\item \label{item:star3} 
Along each solution $\phi$ to $\Sigma_\epsilon$ with $\phi(\dom \phi) \subset U( \partial K) \backslash K$, the map $t \mapsto B(\phi(t))$ is nonincreasing. 
\end{enumerate}

In turn,   \ref{item:star3} holds if there exists a continuous function $\epsilon : \mathbb{R}^n \rightarrow \mathbb{R}_{>0}$ such that  one of the following holds:

\begin{enumerate} [label={C11\alph*.},leftmargin=*]
\item  \label{item:cond1c} $B$ is locally Lipschitz and $F$ is continuous,  and 
$ \langle \nabla B(x) , \eta  \rangle \leq 0$ for all $\eta \in F(x) + \epsilon(x) \mathbb{B}$ and for all $x \in (U(\partial K) \backslash K ) \backslash \Omega$, where $\Omega \subset \mathbb{R}^n$ is the set of points where $B$ is not differentiable.

\item \label{item:cond1b}  $B$ is locally Lipschitz and regular\footnote{\blue{$B$  is \textit{regular}  if 
$  \displaystyle{\lim_{h \rightarrow 0^+}}  \frac{ B(x + hv)  - B(x)}{h} = \max \{ \langle v,  \zeta  \rangle  : \zeta \in \partial_C B(x)  \}.  $}},  and 
\begin{align*}
\hspace{-0.6cm} \langle \partial_C B(x), \eta \rangle \subset  \mathbb{R}_{\leq 0}  \quad \forall \eta \in F_{\epsilon B}(x), 
~  \forall x \in  U(\partial K) \backslash K,
\end{align*}
where $ F_{\epsilon B} :  \mathbb{R}^n  \rightrightarrows \mathbb{R}^n$ is given by
\begin{align*}
\hspace{-1.2cm} F_{\epsilon B}(x)  :=  \{ \eta \in F(x) + \epsilon(x) \mathbb{B} : \text{Card} \left(\langle \partial_C B(x), \eta \rangle \right) = 1   \}.
\end{align*}

\item \label{item:cond1a}  The map $F$ is locally Lipschitz and 
\begin{align*} 
\hspace{-0.6cm} \langle \zeta , \eta  \rangle  \leq 0  \qquad  &  \forall \zeta \in \partial_P B(x),~\forall \eta \in F(x) + \epsilon(x) \mathbb{B}, 
 \\ & \forall x \in  U(\partial K) \backslash K.
\end{align*}
\end{enumerate}}
\end{proposition}

\blue{ It is important to note that  conditions  \ref{item:cond1c}-\ref{item:cond1a} involve the perturbation term $\epsilon$.   This is different from some classical results on robustness analysis \cite{COCV_2000__5__313_0,   goebel2012hybrid}, where the proposed conditions involve only the map $F$. Note that to verify one  of the conditions in \ref{item:cond1c}-\ref{item:cond1a},  when the perturbation is unknown,  we need  to check the inequalities therein for different values of $\epsilon$ until valid candidate is found.  Although being more tedious to check,  such conditions however can lead to larger values of the robustness margin compared to robustness-margin-free conditions that we propose next.   }

In the following result, we establish a set of sufficient conditions for robust safety using barrier functions.

\begin{theorem} \label{thm1}
Consider system $\Sigma$ in \eqref{eq.1} such that Assumption \ref{ass1} holds.
Let $B : \mathbb{R}^n \rightarrow \mathbb{R}$ be a barrier function candidate with respect to $(X_o, X_u) \subset \mathbb{R}^n \times \mathbb{R}^n$. Then, $\Sigma$ is robustly safe with respect to $(X_o,X_u)$ if one the following conditions holds:
\begin{enumerate} [label={C\ref{thm1}\arabic*.},leftmargin=*]
\item \label{item:cond1} There exists  $\epsilon : \mathbb{R}^n \rightarrow \mathbb{R}_{>0}$ continuous such that the set $K$ in \eqref{eq.4} is forward invariant for $\Sigma_{\epsilon}$.

\item \label{item:cond3} The function $B$ is locally Lipschitz and 
\begin{align}
\hspace{-0.8cm} \langle \zeta, \eta \rangle < 0 \quad \forall \eta \in F(x), ~ \zeta \in \partial_C B(x), ~  \forall x \in  \partial K.
\end{align} 

\item \label{item:cond2} The function $B$ is continuously differentiable and 
\begin{align} \label{eq.scond}
\langle \nabla B (x),  \eta \rangle < 0 \quad \forall \eta \in F(x), \quad  \forall x \in  \partial K.
\end{align} 
\end{enumerate}
\end{theorem}

\begin{proof}
Under \ref{item:cond1}, we conclude the existence of a continuous function $\epsilon : \mathbb{R}^n \rightarrow \mathbb{R}_{>0}$ such that the solutions to $\Sigma_{\epsilon}$ starting from $X_o$ never reach $X_u$, which by definition implies robust safety of $\Sigma$ with respect to $(X_o,X_u)$. 

Under \ref{item:cond3}, we introduce the set-valued map $G : \partial K \rightrightarrows \mathbb{R}_{<0}$ given by 
\begin{align} \label{eqGlip}
G(x) := \{ \langle \zeta,  \eta \rangle : \eta \in F(x),~ \zeta \in \partial_C B (x) \}. 
\end{align}
Since $\partial_C B : \mathbb{R}^n \rightrightarrows 
\mathbb{R}^n $ is upper semicontinuous and $\partial_C B(x)$ is nonempty, compact, and convex for all $x \in \mathbb{R}^n$ and since $F$ satisfies Assumption \ref{ass1}, it follows  that 
$G$ satisfies the same properties as $F$. Hence, using Lemma \ref{lemA3} in the Appendix, we conclude the existence of a continuous function $g : \partial K \rightarrow \mathbb{R}_{>0}$ such that 
\begin{align}
\label{eqn:thm11}
\hspace{-0.4cm} \langle \zeta, \eta \rangle \leq  -g(x) ~ \forall \eta \in F(x), ~ \forall \zeta \in \partial_C B(x), ~ \forall x \in  \partial K.
\end{align}
Next, we introduce the function $l : \partial K \rightarrow \mathbb{R}_{\geq 0}$ given by 
\begin{equation}
\label{eqn:l}
 l(x) := \sup \{ |\zeta| : \zeta \in 
\partial_C B(x) \}
\end{equation}
and the function $h: \partial K \rightarrow \mathbb{R}_{>0}$ given by
$ h(x) := \frac{g(x)}{2\sqrt{n}l(x)}. $
The map $\partial_C B$ is upper semicontinuous with compact images, which implies, using \cite[Theorem 1.4.16]{aubin2009set}, that $l$ is upper semicontinuous. Furthermore, under \eqref{eqn:thm11}, we conclude that $l$ is nowhere-zero on $\partial K$. Hence, $h$ is well defined, lower semicontinuous, $-h$ is upper semicontinuous, and $-h(x) <0$ for all $x \in \partial K$. 
Now, using Lemma~\ref{lemA4} in the Appendix, we conclude the existence of a continuous function $\epsilon: \partial K \rightarrow \mathbb{R}_{>0}$ such that $\epsilon(x) \leq h(x)$ for all $x \in \partial K$. Let us show that $\epsilon: \partial K \rightarrow \mathbb{R}_{>0}$ is a robustness margin.  First, we show that
\begin{align*}
\langle \zeta, \epsilon(x)\mu \rangle \leq \frac{g(x)}{2} < 0 \quad  \forall \mu \in \mathbb{B},~ \forall \zeta \in \partial_C B(x), ~ \forall x \in  \partial K,
\end{align*} 
{by noting that} 
 \begin{align*}
 & \langle \zeta, \epsilon(x) \mu \rangle
  = \sum\limits_{i=1}^n \varepsilon(x) \zeta_i \mu_i
 \leq  \epsilon(x) \sum\limits_{i=1}^n |\zeta_i \mu_i|
 \\ & \leq \epsilon(x) \sum\limits_{i=1}^n |\zeta_i| \leq \epsilon(x)  \sqrt{n}|\zeta|
\leq \frac{g(x)}{2\sqrt{n}l(x)}\sqrt{n}|\zeta| \leq  \frac{g(x)}{2},
 \end{align*}
 where the second inequality follows from the fact that $\mu \in \mathbb{B}$, the third inequality comes from the fact that $|\zeta|_1 \leq \sqrt{n} |\zeta|$ {and the last inequality follows from the definition of the map $l: \partial K \rightarrow \mathbb{R}_{\geq 0}$ in (\ref{eqn:l})}.
Finally,  using (\ref{eqn:thm11}),  we conclude that  
\begin{align*}
    \langle \zeta, \eta \rangle \leq
 -\frac{g(x)}{2} < 0 \quad & \forall \eta \in F(x) + \epsilon(x) \mathbb{B}, \\ & \forall \zeta \in \partial_C B(x), ~ \forall x \in  \partial K,
\end{align*}
which implies forward invariance of the set $K$ for $\Sigma_\epsilon$ using 
\cite[Theorem 6]{draftautomatica}. Hence, robust safety of $\Sigma$ with respect to $(X_o,X_u)$ follows using \ref{item:cond1}.  

\blue{Under \ref{item:cond2},  the proof follows {directly} since \ref{item:cond2} implies \ref{item:cond3}.}
\end{proof}

{We now revisit Example \ref{sec.3-a} and show a way to construct a 
triggering strategy that ensures  safety,  using the result of Theorem \ref{thm1}.}

\begin{example}[Safety via self-triggered control] \label{expSTR}
\blue{Let us assume the existence of a (smooth) barrier function candidate $B$ such that \ref{item:cond2} holds. 
A consequence of which, we can find smooth functions $\alpha : \mathbb{R}^n \rightarrow \mathbb{R}$, with $\alpha(x) > 0$ for all $x \in \partial K$, and $\gamma_1 : \mathbb{R}^n \times \mathbb{R}^n \rightarrow \mathbb{R}$, with $\gamma_1 (x,x) = 0$,  such that, for each $(x,y) \in K \times K$,
\begin{equation}
\label{eqrobineq}
\begin{aligned} 
\langle \nabla B(x),  f(x,\kappa(y)) \rangle &  \leq - \alpha(x) + \gamma_1(x,y).
\end{aligned}    
\end{equation} 
Assume,  additionally, that we can find $T_1 > 0$ and $\beta > 0$ such that, for each solution $\phi$ to $\dot{x} = f(x,\kappa(y))$ starting from 
$ y \in \{x \in K: |x|_{\partial K} \geq \beta\}$
satisfies $\phi([0,T_1])  \subset K$. 
Then, using \cite[Theorem 2]{9483096}, we conclude that the problem is solved by simply taking $t_{i+1} - t_i = \min \{T_1,T_r(\phi(t_i))\}$, where, for $y \in K$, we have 
\\
$T_r(y) := \left\{ 
\begin{matrix}
1 & \hspace{-0.4cm} \text{if} ~ M_r(1,y) \leq 0
\\ 
\min \left\{ 1, \frac{2\alpha(y)}{M_r(1,y)} \right\} & \mbox{otherwise},
\end{matrix}
\right.$
\begin{equation*}
\begin{aligned}
M_r(1, y) & := 
\\ &  \sup \{\langle \nabla_x \gamma_1(x, y), f(x, \kappa(y)) \rangle:  x \in Reach(1, y)\}
\\ & + \sup \{\langle - \nabla \alpha(x) , f(x, \kappa(y)) \rangle: x \in Reach(1, y)\},
\end{aligned}
\end{equation*}
and $Reach$ is an over-approximated reachable set  from $y$ over $[0,1]$ along the solutions to $\dot{x} = f(x,\kappa(y))$.}

Finally, Zeno behaviors are excluded provided that
\begin{align}  \label{ineqkey}
T_2 := \inf  \left\{T_r(y)  : y \in K,~|y|_{\partial K} 
\leq \beta \right\} > 0. 
\end{align}
\end{example}
\begin{remark} 
\blue{When a control system  admits a strict control barrier function, as defined in \cite{usevitch2020strong}, then there exists a controller so that the closed-loop system satisfies  \ref{item:cond2}. Indeed, \ref{item:cond2} induces a property named \textit{strict invariance}, or \textit{contractivity}, of the set $K$; namely, in addition to $K$ being forward invariant, the solution  cannot remain in $\partial K$ for a non-trivial time interval.  }
\end{remark}

\begin{remark}
When $B$ is only lower semicontinuous and $F$ is locally Lipschitz,  one can attempt to characterize robust safety  using a condition of the form  
\begin{align}  \label{eqlscadddd}
\langle \zeta, -\eta \rangle < 0  \quad 
\forall (\eta,\zeta,x) \in F(x) \times \partial_P B (x) \times  \partial K.
\end{align}  
Note \eqref{eqlscadddd} allows to conclude robust 
safety provided that,  additionally,  $B$ is convex.  In this case,  
we know that $\partial_P B$ enjoys the same regularity properties as $\partial_C B$; namely,  $\partial_P B$ is upper semicontinuous, with non empty and compact images \cite{cortes2008discontinuous}.  Thus,  the arguments used to prove robust safety under \ref{item:cond3} can be applied.   However,  when $B$ is not convex,  the absence of regularity of 
 $\partial_P B$ prevents the construction of a continuous robustness margin.     
\end{remark}

\begin{remark}
Note that one could 
think of making conditions \ref{item:cond1c}-\ref{item:cond1a} independent of $\epsilon$ by proceeding analogously to \ref{item:cond3} and \ref{item:cond2}.  Namely,  by replacing the $``\leq 0"$  in \ref{item:cond1c}-\ref{item:cond1a} by $``< 0"$ and checking only for vector fields in $F$.  Unfortunately,  when doing so, we fail to guarantee robust safety since the map in\eqref{eqGlip} is not guaranteed to be strictly negative on $\partial K$ and can fail to be upper semicontinuous. 
\end{remark}

In the following result, we propose conditions to guarantee uniform robust safety.

\begin{theorem} \label{thm2}
Consider system $\Sigma$ in \eqref{eq.1} such that Assumption \ref{ass1} holds. 
Let $B : \mathbb{R}^n \rightarrow \mathbb{R}$ be a barrier function candidate with respect to $(X_o, X_u) \subset \mathbb{R}^n \times \mathbb{R}^n$. Then, $\Sigma$ is uniformly robustly safe with respect to $(X_o,X_u)$ if one the following conditions holds:

\begin{enumerate} [label={C\ref{thm2}\arabic*.},leftmargin=*]

\item \label{item:cond1u} There exists a constant $\epsilon > 0$ such that the set $K$ in \eqref{eq.4} is forward invariant for $\Sigma_{\epsilon}$.

\item \label{item:cond2u} 
\blue{The set $\partial K$ is bounded and there exists a continuous function $\epsilon : \mathbb{R}^n \rightarrow \mathbb{R}_{>0}$ such that,  for  $\Sigma_\epsilon$,  the set $K$ is forward invariant.}.

\item \label{item:cond2ubis} 
$\Sigma$ is robustly safe with respect to $(X_o,X_u)$ and either the complement of $X_u$ is bounded or the complement of $X_o$ is bounded.

\item \label{item:cond3u} The function $B$ is continuously differentiable and  
\begin{align*}
\hspace{-0.4cm} \inf \left\{  \frac{\langle \nabla B (x), -\eta \rangle}{|\nabla B (x)|} : (\eta, x) \in  F(x) 
 \times \partial K \right\} > 0.
\end{align*}

\item \label{item:cond4u} The function $B$ is locally Lipschitz and  
\begin{align*}
\hspace{-0.4cm} \inf \left\{ \frac{\langle \zeta, -\eta \rangle}{|\zeta|} : 
 (\eta,\zeta,x) \in F(x) \times \partial_C B(x)  \times \partial K \right\} > 0.
\end{align*} 

\item \label{item:cond5u} 
 $B$ is lower semicontinuous, $F$ is locally Lipschitz, and  
\begin{align*}
\hspace{-1cm} \inf \left\{ \frac{\langle \zeta, -\eta \rangle}{|\zeta|}  :
 (\eta,\zeta,x) \in F(x) \times \partial_P B (x) \times  U(\partial K) \right\}  > 0.
\end{align*}

\item \label{item:cond6u} $B$ is upper semicontinuous, $F$ is locally Lipschitz, $\cl(K) \cap X_u = \emptyset$, and  
\begin{align*}
\hspace{-1cm} \inf  \left\{ \frac{\langle \zeta, \eta \rangle}{|\zeta|}  : 
 (\eta,\zeta,x) \in F(x) \times \partial_P (-B(x)) \times  U(\partial K) \right\}  > 0.
\end{align*}
\end{enumerate}
\end{theorem}

\begin{proof}
Under \ref{item:cond1u}, we conclude the existence of $\epsilon > 0$ such that the solutions to $\Sigma_{\epsilon}$ starting from $X_o$ never reach $X_u$, which by definition implies uniform robust safety of $\Sigma$ with respect to $(X_o,X_u)$.

\blue{ Under \ref{item:cond2u} and when $\partial K$ is bounded,  we show that 
$\epsilon^* := \inf \{ \epsilon(x) : x \in U(\partial K) \}$, for $U(\partial K)$ a neighborhood of $\partial K$,   is a robustness margin using contradiction.  Indeed,  assume that there exists a solution 
$\phi$ to $\Sigma_{\epsilon^*}$ starting from $x_o \in X_o$ that reaches the set $X_u$ in finite time.  
Namely,  there exists $t_1$,  $t_2 \in \dom \phi$ such that $t_2 > t_1 \geq 0$,  $\phi(t_2) \in U(\partial K) \backslash K$,  $\phi(t_1) \in U(\partial K) \cap K$, 
$\phi(0) \in X_o$,  and $\phi([t_1,t_2]) \subset U(\partial K)$.  The latter implies that $\phi$,  restricted to the interval $[t_1,t_2]$,  is also a solution to 
$\Sigma_\epsilon$.   However, since the set $K$ is forward invariant for $\Sigma_{\epsilon}$,  we conclude that the solution $\phi([t_1,t_2])$ must lie within the set $K$, which yields to a contradiction.   }

Under \ref{item:cond2ubis}, we conclude the existence of $\epsilon : \mathbb{R}^n \rightarrow \mathbb{R}_{>0}$ continuous such that $\Sigma_\epsilon$ is safe with 
respect to $(X_o,X_u)$. Hence, the set 
$ K_{\epsilon} := \{ \phi(t) : t \in \dom \phi,~\phi \in \mathcal{S}_{\Sigma_\epsilon}(x),~ x \in X_o \} $ 
is forward invariant for $\Sigma_\epsilon$, $K_{\epsilon} \cap X_u = \emptyset$, and $X_o \subset K_{\epsilon}$. Furthermore, when either the complement of $X_u$ is bounded or the complement of $X_o$ is bounded., we conclude that $\partial K_{\epsilon}$ is bounded. Hence, the rest of the proof follows as in the proof of \ref{item:cond2u}.

Under \ref{item:cond3u},  
 \blue{we choose 
$$ \epsilon := \inf \left\{ \frac{\langle \nabla B (x), -\eta \rangle}{2\sqrt{n} |\nabla B (x)| } : (\eta, x) \in  F(x) \times \partial K\right\}, $$
Consider $x \in \partial K$ and $\mu \in \mathbb{B}$, we have:
\begin{align}
& \nonumber \langle \nabla B(x),  \epsilon\mu \rangle  = \sum\limits_{i=1}^n\epsilon\nabla_{x_i} B(x) \mu_i 
\\ \label{eqn:proof1} & \leq \epsilon \sum\limits_{i=1}^n 
 |\nabla_{x_i} B(x)|    \leq \epsilon  \sqrt{n}|\nabla B(x)|,
 \end{align}
 where the second inequality follows from the fact that $\mu \in \mathbb{B}$ and the third inequality comes from the fact that 
 $ |x|_1 := \sum_{i=1}^n|x_i| \leq \sqrt{n} |x|$. Now consider $x \in \partial K$, $\nu \in F(x)$ and $\mu \in \mathbb{B}$ that 
\begin{align*}
\langle \nabla B (x),&  \nu+\epsilon \mu \rangle = \langle \nabla B (x),  \nu\rangle +\langle \nabla B (x),  \epsilon \mu \rangle \\& \leq \langle \nabla B (x),  \nu\rangle +  \epsilon  \sqrt{n}|\nabla B(x)| \\ &\leq \langle \nabla B (x),  \nu\rangle + \frac{\langle \nabla B (x), -\nu \rangle}{2\sqrt{n} |\nabla B (x)}|\sqrt{n}|\nabla B(x)| \\& \langle \nabla B (x),  \frac{\nu}{2}\rangle < 0, 
\end{align*}
where the second inequality comes from (\ref{eqn:proof1}) and the last inequality follows from the definition of $\epsilon$. Hence, we conclude that
\begin{align*}
\langle \nabla B (x),  \eta \rangle < 0 
\qquad \forall \eta \in F(x) + \epsilon \mathbb{B},  \quad \forall x \in  \partial K,
\end{align*}}
which, in turn, implies forward invariance of the set $K$ for $\Sigma_\epsilon$ using {\cite[Theorem 5]{draftautomatica}}. Hence, uniform robust safety of $\Sigma$ with respect to $(X_o,X_u)$ follows using \ref{item:cond1u}.

Similarly, under \ref{item:cond4u} and for 
$$ \epsilon := \inf \left\{
\frac{\langle \zeta, -\eta \rangle}{\sqrt{n} |\zeta|} : (\eta,\zeta, x) \in  F(x) \times \partial_C B(x) \times \partial K \right\}, $$
we conclude that 
\begin{align*}
\langle \zeta,  \eta \rangle < 0 
\quad \forall \eta \in F(x) + \epsilon \mathbb{B}, ~\forall \zeta \in \partial_C B(x), ~  \forall x \in  \partial K,
\end{align*} 
which, in turn, implies forward invariance of the set $K$ for $\Sigma_\epsilon$ using \cite[Theorem 6]{draftautomatica}. Hence, uniform robust safety of $\Sigma$ with respect to $(X_o,X_u)$ follows using \ref{item:cond1u}.

Now, under \ref{item:cond5u} and for 
$$ \epsilon := \inf \left\{
\frac{\langle \zeta, -\eta \rangle}{\sqrt{n} |\zeta|} : (\eta,\zeta, x) \in  F(x) \times \partial_P B(x) \times U(\partial K) \right\}, $$
we conclude that 
\begin{align} \label{eqdecrease}
\langle \zeta,  \eta \rangle \leq 0 
\quad \forall \eta \in F(x) + \epsilon \mathbb{B}, ~\forall \zeta \in \partial_P B(x), ~  \forall x \in  U(\partial K).
\end{align} 
Now, using Lemma \ref{lemA9}, we conclude that \eqref{eqdecrease} implies that, along each solution $\phi$ to $\Sigma_\epsilon$ starting from $x_o \in U(\partial K)$ and remaining in $U(\partial K)$, the map $t \mapsto B(\phi(t,x_o))$ is nonincreasing. 
The latter implies, using Lemma \ref{lemA8}, that $\Sigma_\epsilon$ is safe with respect $(X_o,X_u)$. Hence, uniform robust safety of $\Sigma$ with respect to $(X_o,X_u)$ follows using \ref{item:cond1u}.

Finally,  under \ref{item:cond6u}, we take 
$$ \epsilon := \inf \left\{
\frac{\langle \zeta, -\eta \rangle}{\sqrt{n} |\zeta|} : (\eta,\zeta, x) \in  F(x) \times \partial_P B(x) \times U(\partial K) \right\}, $$
to obtain 
\begin{equation}
\label{eqdecrease1}
\begin{aligned} 
\langle \zeta,  \eta \rangle \leq 0 
\quad \forall \eta \in -F(x) + \epsilon \mathbb{B}, ~ & \forall \zeta \in 
\partial_P(-B(x)), \\ &  \forall x \in  U(\partial K).
\end{aligned} 
\end{equation}
Using Lemma \ref{lemA9}, we conclude that \eqref{eqdecrease1} implies that, along each solution $\phi$ to $\Sigma^-_\epsilon$ starting from $x_o \in U(\partial K)$ and remaining in $U(\partial K)$, the map $t \mapsto 
-B(\phi(t,x_o))$ is nonincreasing, where 
$ \Sigma^{-}_{\epsilon} : \quad \dot{x} \in - F(x) + \epsilon \mathbb{B}$.  
Hence, along each solution $\phi$ to $\Sigma_\epsilon$ starting from $x_o \in U(\partial K)$ and remaining in $U(\partial K)$, the map $t \mapsto B(\phi(t,x_o))$ is nonincreasing. The latter implies, using Lemma \ref{lemA8}, that $\Sigma_\epsilon$ is safe with respect $(X_o,X_u)$. Hence, uniform robust safety of $\Sigma$ with respect to $(X_o,X_u)$ follows using \ref{item:cond1u}. 
\end{proof}

\begin{remark}
The sufficient conditions for uniform robust safety in Theorem \ref{thm2} generalize those proposed in \cite{ratschan2018converse} and \cite{wisniewski2016converse}, where only \ref{item:cond2ubis} is considered in the particular case of differential equations.
\end{remark}

\section{Numerical Examples} \label{sec.exp}

\ifitsdraft
\red{
In the following example,  we propose a system that is safe but not robustly safe.   

\begin{example} 
\label{exp1} Consider the system
\begin{equation*}
\begin{aligned}
\Sigma : 
\left\{
\dot{x} =  F(x) = 
\begin{bmatrix}
0 & 1
\\
-1 & 0
\end{bmatrix}
x
\qquad x \in \mathbb{R}^2.
\right.
\end{aligned}
\end{equation*}
The initial and unsafe sets are given by 
\begin{align} \label{eqsets}
X_o := \{x \in \mathbb{R}^2 : x_1^2+x_2^2 \leq 1 \} ~ \text{and} ~
X_u :=  \mathbb{R}^2 \backslash X_o.
\end{align}
Suppose that $\Sigma$ is affected by a disturbance leading to
\begin{equation} \label{eqsigmaeps}
\begin{aligned}
\Sigma_{\epsilon} : 
\left\{
\dot{x} = 
\begin{bmatrix}
0 & 1
\\
-1 & 0
\end{bmatrix}
x + 
\begin{bmatrix}
\epsilon 
\\
0
\end{bmatrix}
\qquad x \in \mathbb{R}^2,
\right.
\end{aligned}
\end{equation}
where $\epsilon \geq 0$ is constant. First, we note that, when $\epsilon = 0$,  system $\Sigma$ is safe with respect to $(X_o,X_u)$. Indeed, using the barrier function candidate $ B(x) := x_1^2+x_2^2 - 1, $ we conclude that 
$\langle \nabla B (x), F(x) \rangle \leq 0$ for all $x \in \mathbb{R}^2$. 
Hence, safety follows 
using Lemmas \ref{lemA8} and \ref{lemA9}. However, as soon as $\epsilon$ becomes strictly positive, safety is violated; see Figure \ref{fig:example5}. Indeed, the solution to $\Sigma_\epsilon$ starting from $x_o := (x_{o1},x_{o2})$, when $\epsilon$ is constant, is given by
\begin{equation*}
\begin{cases}
\phi_1(t)  = & (x_{o2} + \epsilon) \sin(t) + x_{o1} \cos(t) 
\\ 
\phi_2(t)  = & (x_{o2} 
+ \epsilon) \cos(t) -  x_{o1} \sin(t) -  \epsilon.
\end{cases} 
\end{equation*}
Hence, for $x_o:=(0,1) \in X_o$, we have 
\begin{equation*}
\begin{cases}
\phi_1(t)  = & (1 + \epsilon) \sin(t) 
\\ 
\phi_2(t)  = & (1 
+ \epsilon) \cos(t) -  \epsilon
\end{cases} 
\end{equation*}
with $\phi(\pi/2) = (1+\epsilon, -\epsilon) \in X_u$. 
\begin{figure}[t]
\centering
  \includegraphics[scale=0.21]{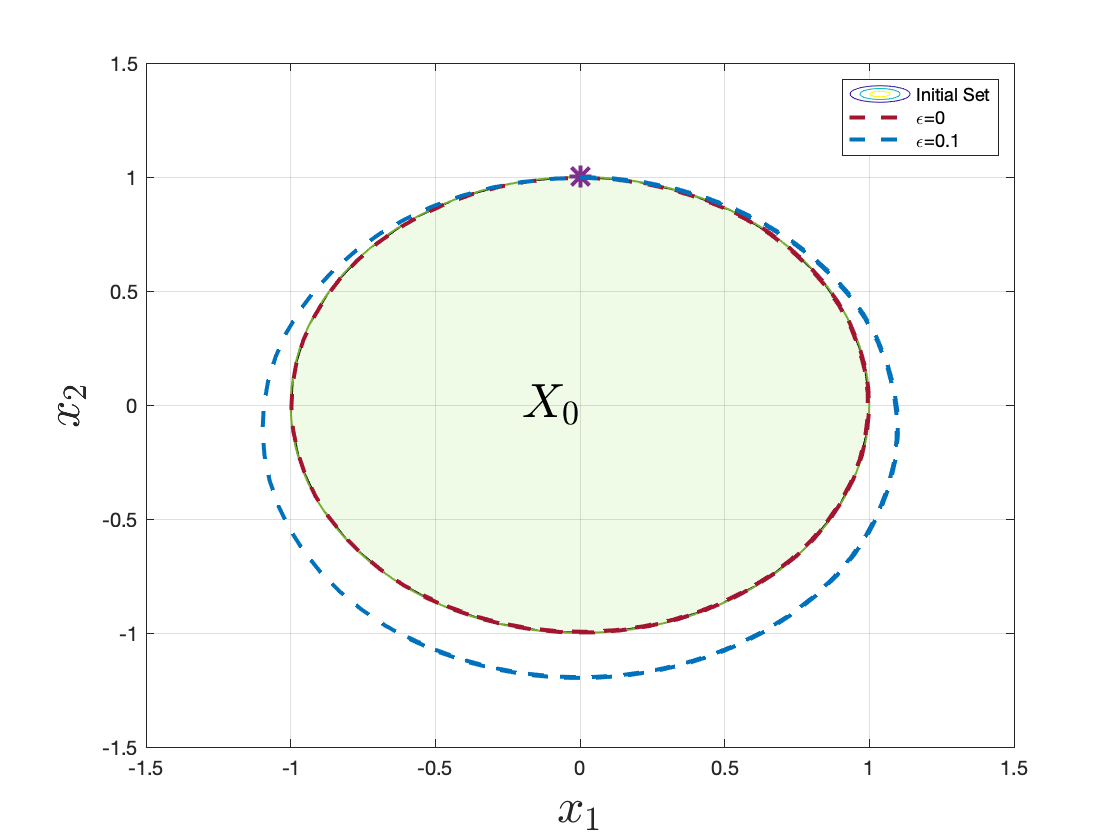}
  \caption{\textcolor{red}{Illustration of the sets $X_o$ and $X_u$ and trajectories of $\Sigma$ and $\Sigma_{\epsilon}$, for different values of $\epsilon$, initiated from $x_o := (0,1)$. } }
  \label{fig:example5}
\end{figure}
\end{example}

\fi

\begin{example}[Safety via 
self-triggered control] \label{exp...}
\blue{Consider the control system 
$$ \Sigma_u : \dot{x} = \begin{bmatrix}  u & & - x_1 \end{bmatrix}^\top \qquad x \in \mathbb{R}^2, \quad u \in \mathbb{R}.  $$
Consider the feedback law  $\kappa_o(x) := x_2$ rendering the closed-loop system 
$ \Sigma : \dot{x} = f(x, \kappa_o(x)) = \begin{bmatrix}  x_2 & & - x_1 \end{bmatrix}^\top  $
safe with respect to 
\begin{align} \label{eqsets}
X_o := \{x \in \mathbb{R}^2 : x_1^2+x_2^2 \leq 1 \} ~ \text{and} ~
X_u :=  \mathbb{R}^2 \backslash X_o.
\end{align}
Given a triggering sequence $\{ t_i \}^\infty_{i=1}$ of period $T > 0$, a solution $\phi$ to the self-triggered closed-loop system is governed, for all $t \in [t_i, t_{i+1})$, by 
\begin{align} \label{eqsTcL}
\hspace{-0.4cm} 
\dot{\phi}(t) = 
\begin{bmatrix} 
0 & 1 \\ -1 & 0
\end{bmatrix} \phi(t) + 
\begin{bmatrix} 
\phi_2(t_i) - \phi_2(t) \\ 0
\end{bmatrix}. 
\end{align}
Since we already established in \cite[Example 1]{9683684} that 
\begin{equation*}
\begin{aligned}
\Sigma_{\epsilon} : 
\left\{
\dot{x} = 
\begin{bmatrix}
0 & 1
\\
-1 & 0
\end{bmatrix}
x + 
\begin{bmatrix}
\epsilon 
\\
0
\end{bmatrix}
\qquad x \in \mathbb{R}^2,
\right.
\end{aligned}
\end{equation*} 
cannot be safe even if $\epsilon >0$ is arbitrary small, one should expect that, for any $T > 0$ {arbitrarily} small, there exists a solution $\phi$ to the self-triggered closed-loop system, starting from $X_o$, that leaves the set $X_o$; see Figure \ref{fig:example6}. Indeed, let $\phi(0) := (0,1)$, solving \eqref{eqsTcL}, we conclude that, for all $t \in [0,T]$, $|\phi(t)|^2 = 1 + \frac{t^4}{4} > 1$.  }
\begin{figure}[t]
\centering
\includegraphics[scale=0.18]{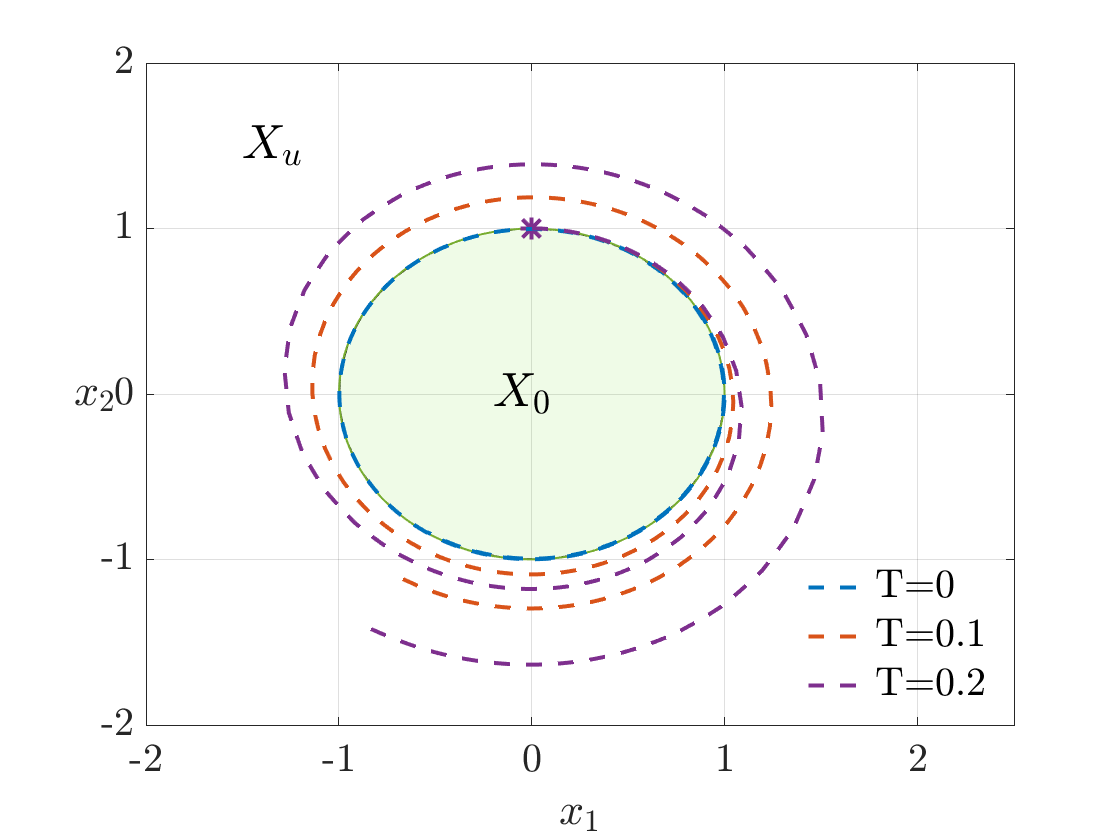}
  \caption{  Illustration of $X_o$ and $X_u$, a  trajectory to $\Sigma_u$ subject to a continuous implementation of $\kappa_o$, and trajectories of $\Sigma_u$ subject to a self-triggered implementation of $\kappa_o$, for different values of $T > 0$.}
  \label{fig:example6}
\end{figure}

\blue{We, now, use the controller 
$\kappa(x) := - x_1 + x_2$, which renders the closed-loop system $\Sigma : \dot{x} = f(x, \kappa(x)) = \begin{bmatrix}  -x_1 + x_2 & & - x_1 \end{bmatrix}^\top$
robustly safe with respect to 
$ X_o := \{x \in \mathbb{R}^2 : W(x) \leq 1 \}$ and
$X_u :=  \mathbb{R}^2 \backslash X_o$, where 
\begin{align} \label{eqV}
W(x) := x_1^2+x_2^2 - x_1 x_2.
\end{align}
Indeed, it is easy to see that \ref{item:cond2} holds for $B(x) := W(x) - 1$. 
Now, following Example \ref{expSTR}, we show the existence of a periodic triggering sequence that maintains the system's safety.}
\blue{Note that \eqref{eqrobineq} is satisfied for $\alpha(x) := W(x)$ 
and $\gamma_1(x,y) :=  (2 x_1 - x_2)  \left[(x_2 - y_2)-(x_1 - y_1) \right]$. Next, since  
$\sup \{|\gamma_1(x,y)| : W(x) \leq 1, W(y) \leq 1/2 \} \leq 21$, we conclude that, for $y$ satisfying $W(y) = 1/2$ and along the solutions to $\dot{x} = f(x,\kappa(y))$, we have 
$\dot{W} \leq - W + 21$ as along as $W \leq 1$. Hence, we can take $T_1 = 0.02$ and $\beta := \max \{|y|_{\partial K} : W(y) = 1/2 \} \leq 7/4$ to obtain $T_2 = 0.19$. Finally, to show that  \eqref{ineqkey} holds, we can just use 
\cite[Theorem 4]{9483096}, 
which, additionally, requires  from the set  $\{y \in K : |y|_{\partial K} \leq \beta \}$
to be compact and from the map $y \mapsto Reach(1,y)$ to be upper semicontinuous with compact images.} 

\begin{figure}[t]
\centering
  \includegraphics[scale=0.18]{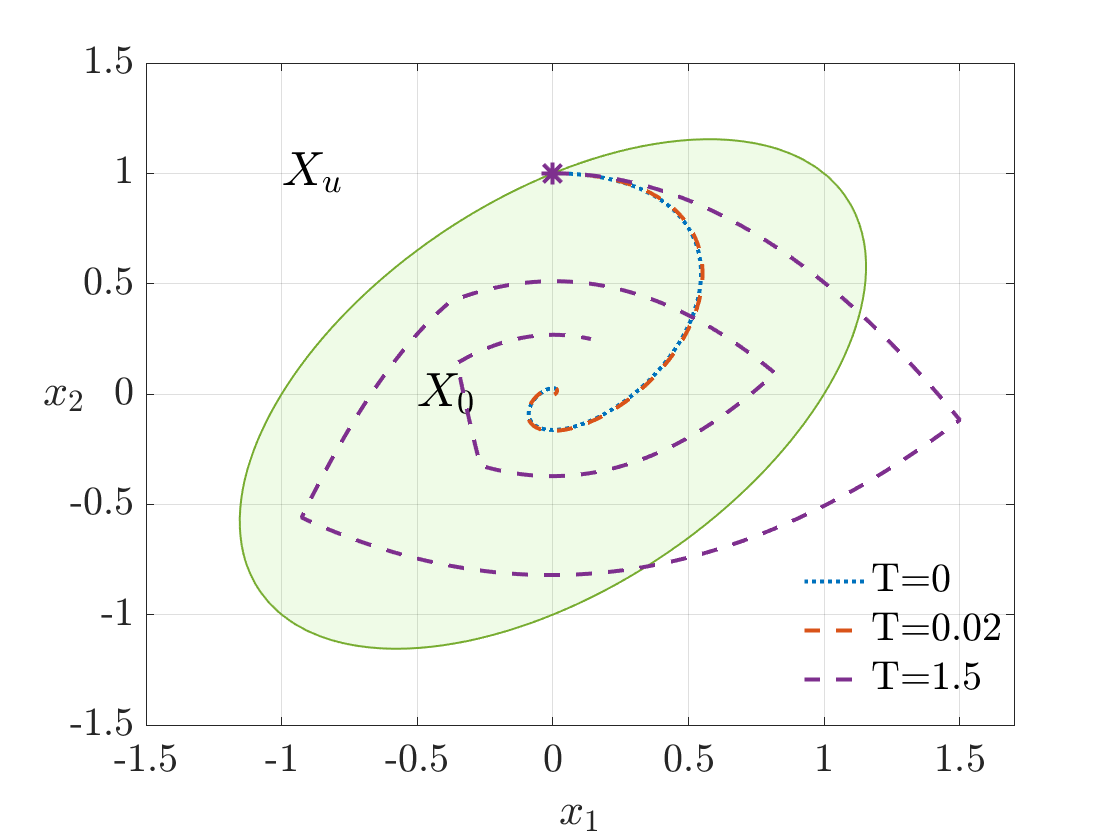}
  \caption{Illustration of $X_o$ and $X_u$, a  trajectory to $\Sigma_u$ subject to continuous implementation of $\kappa$, and trajectories of $\Sigma_u$ subject to self-triggered implementation of $\kappa$ for different values $T$.}
\label{figure:example6bis}
\end{figure}
\end{example}

\begin{example}[Safety under singular perturbations] \label{expsing}
\blue{We revisit Example \ref{expsingpert} by considering the singularly-perturbed system 
 \begin{align} \label{eqsingpert} 
 \dot{x} = 
 \begin{bmatrix}
 -1 & 1 
 \\
 -1 & 0
 \end{bmatrix}
x +
\begin{bmatrix}
1 \\ 1
\end{bmatrix} e,
\quad 
 \varepsilon \dot{e}  = - e + \varepsilon x_2.
 \end{align}
The slow dynamics is governed by
$\dot{x} =  
\begin{bmatrix}
 -1 & 1 
 \\
 -1 & 0
 \end{bmatrix}
x$, which is, as shown in Example \ref{exp...}, robustly safe with respect to $X_o  := \{ x \in \mathbb{R}^2 : W(x) \leq 1 \}$ and $X_u := \{x \in \mathbb{R}^2  : W(x) > 12 \}$, 
where $W$ is defined in \eqref{eqV}. 

We would like here to show the existence of $\varepsilon^* > 0$ such that, for each $\varepsilon \in (0,\varepsilon^*]$, \eqref{eqsingpert} is safe with respect to $(X_o \times \mathbb{B}, X_u \times \mathbb{R})$; see Figure \ref{exp:SPS}.
To do so, we follow the approach suggested in Example \ref{expsingpert}, which consists in finding $R$, $r$, and $T > 0$ so that steps St1. and St2. therein are verified. Note that  
\eqref{eqgamma3} is verified for $\Gamma_3(e)  := \sup \{ |x_2| : W(x) \leq 12 \} = \text{const}$.  Furthermore, using $V(e) := e^2$, we have $\langle \nabla V(e), g(x,e,0) \leq -2 V(e)$;
thus, St2. is verified for any $R > 1$, $T > 0$, $r > 0$. 

Now, to verify St1., having $\Gamma_1(x,e) := |e|$, we need to find $R > 1$, $r$, and $T > 0$ such that 
$ \dot{x} \in F(x) + r \mathbb{B} $ is safe with respect $(Reach(T, X_o), X_u)$.Note that $Reach(T, X_o)$ is the reachable set from 
 $X_o$ over $[0,T]$ along the solutions to
$ \dot{x}  \in  F(x) + R \mathbb{B}$. Furthermore, along the solutions to the latter system, we have
$\dot{W} \leq - W + R |\nabla W|$. Now, by exploiting the quadratic form of $W$, we conclude that $\dot{W} \leq - (1/2) W + 5 R^2$. 
As a result, $R(T,X_o) \subset \{ x \in \mathbb{R}^2 : W(x) \leq 1 + 10 R^2 \}$.  Hence,  by taking $R^2 := 11/10$, we conclude that
$ R(T,X_o) \subset \{ x \in \mathbb{R}^2 : W(x) \leq 12 \}  = \mathbb{R}^2 \backslash X_u. $ 
As a last step, it is enough to  show that 
$ \dot{x} \in F(x)$ is uniformly robustly safe with respect $(\mathbb{R}^2 \backslash X_u, X_u)$. The latter is true using \ref{item:cond3u} and the barrier candidate $B(x) := W(x) - 12$. }
\begin{figure}[t]
\centering
\includegraphics[scale=0.18]{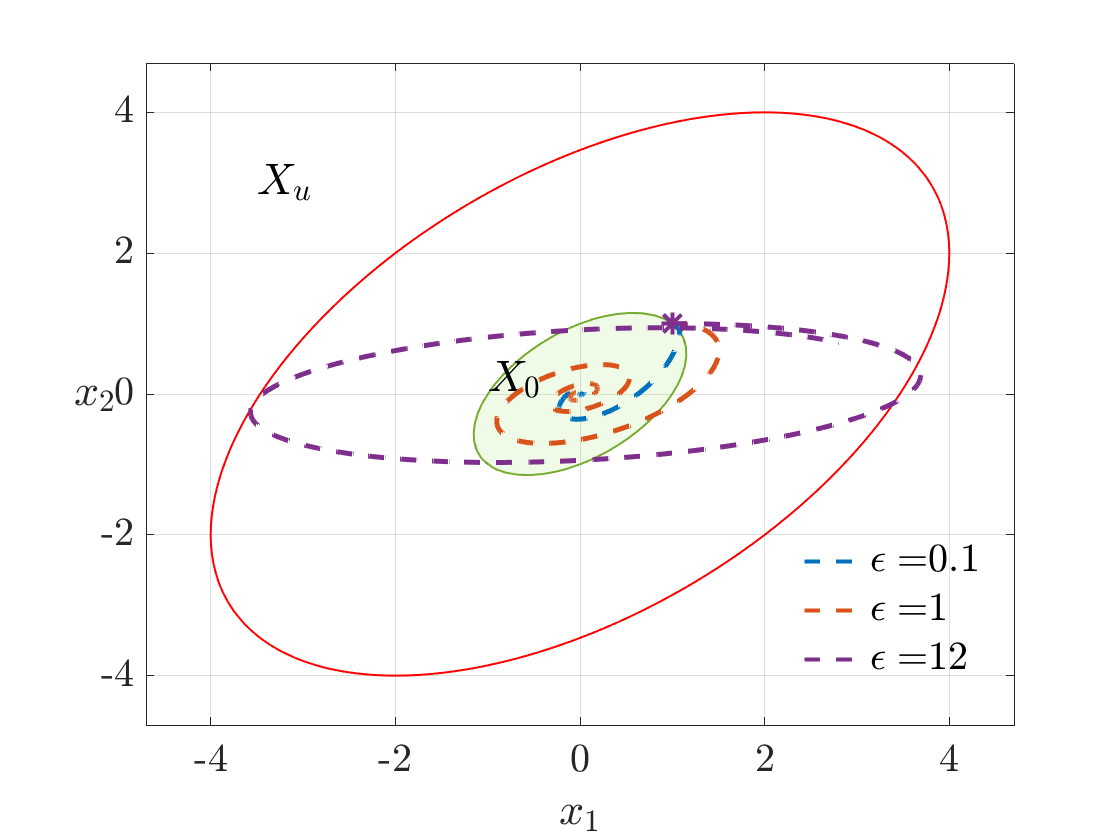}
  \caption{Illustration of $X_o$ and $X_u$ and a  trajectories to the $x$ component of the singularly perturbed system, for different values of $\epsilon$, and for $e(0) = x_1(0) = x_2(0) = 1$. }
  \label{exp:SPS}
\end{figure}

\end{example}

\ifitsdraft

\red{Finally,  we propose a system that is robustly safe but not uniformly robustly safe. 

\begin{example}  \label{exp2}
Consider the system
\begin{equation*}
\begin{aligned}
\Sigma : 
\left\{
\dot{x} =  F(x) = 
\begin{bmatrix}
1
\\
2 x_1 x_2 - \frac{(1 + 2|x_1 x_2|)}{1 + x_1^2 e^{x^2_1}}
\end{bmatrix}
\qquad x \in \mathbb{R}^2.
\right.
\end{aligned}
\end{equation*}

The initial and unsafe sets are given by
\begin{align*}
X_o & := \left\{ x \in \mathbb{R}^2 : x_2 \leq e^{x^2_1} \right\}  \quad \text{and} \quad 
X_u := \mathbb{R}^2 \backslash X_o.
\end{align*}
We will show that system $\Sigma$ admits 
$ \epsilon(x) := \frac{1}{1 + x_1^2 e^{x^2_1}}$
as a robustness margin. Indeed, using the barrier function candidate $B(x) := x_2 e^{-x^2_1} - 1$,  
we conclude that, for each $x \in \mathbb{R}^2$,
$ \langle \nabla B (x), F(x) \rangle = -\epsilon(x) (1+2|x_1x_2|) e^{-x_1^2}$. 
Furthermore,  for $\mu := (\mu_1,\mu_2) \in \mathbb{B}$ and for each $x \in \mathbb{R}^2$, we have
\begin{align*}
 \langle \nabla B(x), \epsilon(x) \mu \rangle & =  \epsilon(x) \left[ \nabla_{x_1} B(x) \mu_1 + \nabla_{x_2} B(x) \mu_2 \right]  \\ &\leq \epsilon(x) (1+2|x_1x_2|)e^{-x_1^2}.
 \end{align*} 
As a result, for each 
$\eta \in F(x) + \epsilon(x) \mathbb{B}$, we conclude that $\langle \nabla B (x), \eta \rangle \leq 0$ for all $x \in \mathbb{R}^2$. Finally, using Lemmas \ref{lemA8} and \ref{lemA9} applied to $\Sigma_\epsilon$ instead of $\Sigma$,  we conclude that $\Sigma_\epsilon$ is safe with respect to $(X_o,X_u)$. Hence, $\Sigma$ is robustly safe with respect to $(X_o,X_u)$ and $\epsilon$ is a robustness margin. 

On the other hand, to conclude that $\Sigma$ is not uniformly robustly safe with respect to $(X_o,X_u)$, we show that, for each constant  
$\epsilon >0$, there exists $x_{\epsilon} \in \partial X_o \cap X_o$ and a solution to $\Sigma_\epsilon$ starting from $x_{\epsilon}$ that reaches $X_u$. Recall that
$$ \langle \nabla B (x), F(x) \rangle = -\frac{1+2|x_1x_2|}{1 + x_1^2 e^{x^2_1}} e^{-x_1^2}. $$
Hence,  for each $x \in \partial X_o \cap X_o$, we have $x_2 = e^{x^2_1}$ and 
\begin{align*} 
& \langle \nabla B (x), F(x) + [0,1]^\top \epsilon \rangle = 
\left(\epsilon -\frac{1 + 
2 |x_1| e^{x_1^2}}{1 + x_1^2 e^{x^2_1}}  \right) e^{-x_1^2}.
\end{align*}
Since the map $x_1 \mapsto \frac{1+|x_1| e^{x_1^2}}{1 +x_1^2 e^{x^2_1}}$ is continuous and converges to $0$ when $x_1$ goes to $\infty$, we conclude that, for any $\epsilon>0$, we can always find 
$ x_\epsilon := (x_{1\epsilon}, x_{2\epsilon}) =(x_{1\epsilon},  e^{x_{1\epsilon}^2}) \in \partial X_o \cap X_o $ 
such that $ \langle \nabla B (x_\epsilon), F(x_\epsilon) + [0,1]^\top \epsilon \rangle > 0$.   Furthermore,  using Lemma   \ref{lemA002},  we conclude that $F(x_\epsilon) + [0,1]^\top \epsilon \in D_{X_u} (x_\epsilon)$,  where $D_{X_u}$ is the \textit{Dubovitsky-Miliutin} cone of $X_u$ defined in \eqref{eq.cone2}.  Finally,  using Lemma \ref{lemA001}, we conclude  the existence of a solution $\phi_\epsilon$, starting from $x_\epsilon$,  to the system  $ \dot{x} = F(x_\epsilon) + [0 \quad 1]^\top \epsilon$, $x \in \mathbb{R}^2$,
which is also a solution to $\Sigma_\epsilon$, that reaches $X_u$.
\begin{figure}[t]
\centering
  \includegraphics[scale=0.21]{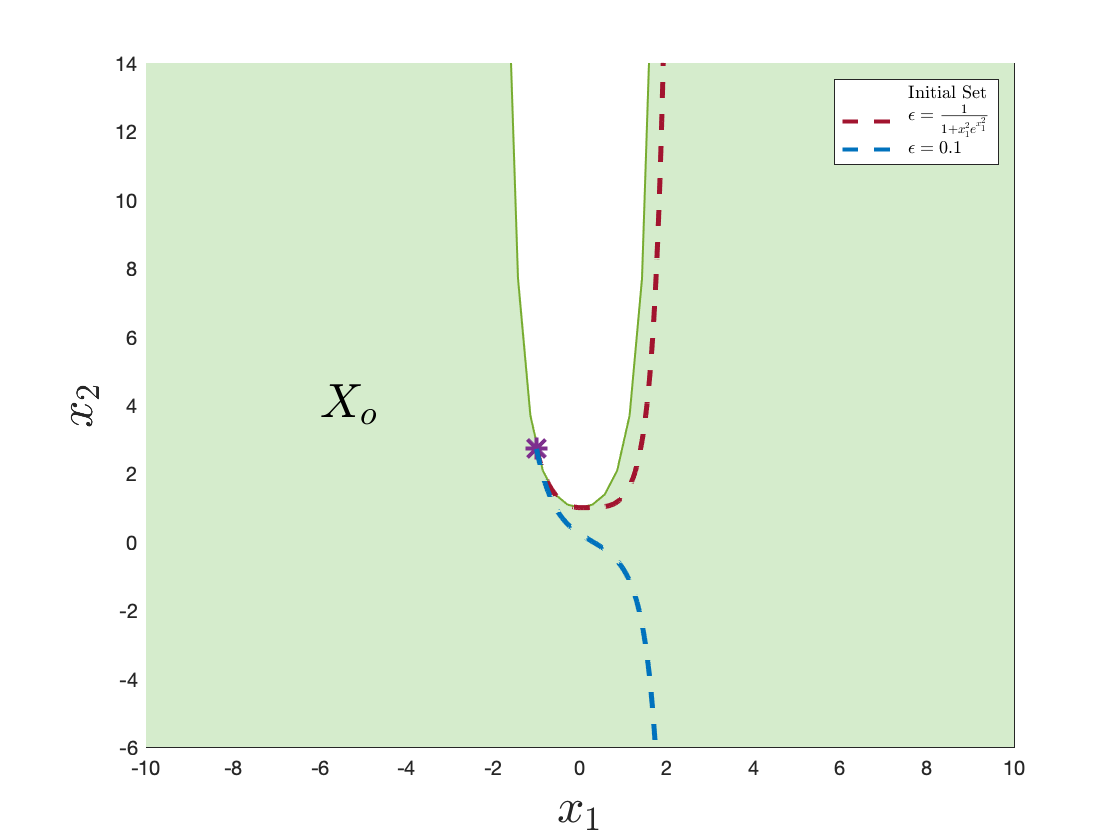}
  \caption{
  \textcolor{red}{ Illustration of $X_o$ and $X_u$, 
  a trajectory to $\Sigma_{\epsilon}$, when $\epsilon \in \{0.1, 0.2\}$, that leaves the set $X_o$, and a trajectory to $\Sigma_{\epsilon(x)}$, $\epsilon(x) := \frac{1}{1 + x_1^2 e^{x^2_1}}$,  that remains within $X_o$.}}
  \label{exp:NUniform}
\end{figure}
\end{example}}

\fi

\begin{example}[Contracts via robust safety] \label{examp:agc}
 \blue{Consider system $\Sigma$, which is the feedback composition of $\Sigma_{w_1}$ and $\Sigma_{w_2}$ given by
\begin{align*} 
\Sigma_{w_1} & : \quad  \dot x_1 = \sqrt{w_1} - 2x_1 -2 \qquad x_1, w_1 \in \mathbb{R}.
\\
\Sigma_{w_2} & : \quad  \dot x_2 = - x_2 + w_2 \qquad ~~ \qquad x_2, w_2 \in \mathbb{R}.
\end{align*}
The feedback composition of  $\Sigma_{w_1}$ 
and $\Sigma_{w_2}$ is constrained by the algebraic equations $w_1=x_2$ and $w_2=x_1$. Hence, $\Sigma$ is given by
\begin{equation*}
\begin{aligned}
\Sigma : \dot{x} = \begin{bmatrix}
\dot x_1
\\
\dot x_2
\end{bmatrix} = 
\begin{bmatrix}
\sqrt{x_2} - 2x_1 -2
\\
-x_2 + x_1
\end{bmatrix}
\qquad x \in \mathbb{R}^2.
\end{aligned}
\end{equation*}
We would like to show that  $\Sigma$ is safe with respect to $(K, \mathbb{R}^2 \backslash K)$, with $K=[0,2] \times [0,2]$. 
To do so, we use \cite[Theorem 2]{saoud2021assume} to conclude that it is enough to show that $\Sigma_{w_2}$ weakly satisfies the contract $\mathcal{C}_2 := (A_2,K_2)$, and $\Sigma_{w_1}$ strongly satisfies the contract $\mathcal{C}_1 := (A_1,K_1)$, where $A_1=A_2=K_1=K_2=[0,2]$. 
First, one can easily show that the component $\Sigma_{w_2}$ weakly satisfies $\mathcal{C}_2$.
To show that the component $\Sigma_{w_1}$ strongly satisfies the contract $\mathcal{C}_1$, we show that the system $\Sigma_1$ defined by
$\Sigma_1 : \dot x_1 \in F_1(x_1):= \sqrt{A_1} -2x_1-2$ is robustly safe with respect to $(K_1,\mathbb{R}^n \backslash K_1)$. Indeed, using the barrier function candidate $B(x_1) := x_1 - 2$,  for each $x_1 \in \partial K_1$, and for each $\eta \in F_1(x_1) $, we have that $\langle \nabla B(x_1), \eta \rangle <
0$. Hence, using the result of Theorem \ref{thm1}, one has that $\Sigma_1$ is robustly safe with respect to $(K_1,\mathbb{R}^n \backslash K_1)$, which implies using Proposition \ref{propagc} that $\Sigma_{w_1}$ strongly satisfies the contract $\mathcal{C}_1=(A_1,K_1)$.}
\end{example}

\section{Conclusion}
The current paper studies two robust-safety notions for differential inclusions. After motivating the considered notions, barrier-function-based conditions are proposed to certify them. Those conditions are infinitesimal as they involve only the unperturbed system's dynamics. In future work, we would like to study the necessity of the proposed sufficient conditions, as well as possible extensions to the general context of hybrid dynamical systems. 
  
\appendix  

\section{Appendix}

\begin{lemma} [Monotonicity implies safety]  \label{lemA8} 
Given initial and unsafe sets
$(X_o, X_u) \subset \mathbb{R}^n \times \mathbb{R}^n$,  system $\Sigma$ is safe with respect to $(X_o, X_u)$ if there exists  a barrier function candidate $B : \mathbb{R}^n \rightarrow \mathbb{R}$ such that 
\begin{enumerate} 
[label={($\star$)},leftmargin=*]
\item \label{item:star1} 
Along each solution $\phi$ to $\Sigma$ with $\phi(\dom \phi) \subset U(\partial K)$,   the map $t \mapsto B(\phi(t))$ is non-increasing. 
\end{enumerate}
Moreover,  if $K$ is closed,  we can replace $U(\partial K)$ in \ref{item:star1} by $U(\partial K) \backslash \text{int}(K)$,  and  if $B$ is continuous,  we can replace $U(\partial K)$ in \ref{item:star1} by $U(\partial K) \backslash K$. 
\end{lemma}

\begin{proof}
\blue{   
To find a contradiction,   we assume the existence of  a solution $\phi$ starting from $x_o \in X_{o}$ that reaches the set $X_{u}$ in finite time,  and we note that the (zero-sublevel) set $K$ satisfies $X_o \subset K$ and $X_u \subset \mathbb{R}^n \backslash K$.  Using continuity of $\phi$,  we conclude the existence of $0 \leq t_1 < t_2$ such that $\phi(t_2) \in  U(\partial K) \backslash K$,  $\phi(t_1) \in K$,  and $\phi([t_1,t_2]) \subset U(\partial K)$.    As a result,  having $B(\phi(t_1)) \leq 0$ and $B(\phi(t_2)) > 0$ contradicts \ref{item:star1}.  
Now,   when $K$ is closed,  we conclude the existence of $0 \leq t_1 < t_2$ such that $\phi(t_2) \in  U(\partial K) \backslash K$,  $\phi(t_1) \in \partial K$,  and $\phi([t_1,t_2]) \subset U(\partial K) \backslash \text{int}(K)$.  
Similarly,  having $B(\phi(t_1)) \leq 0$ and $B(\phi(t_2)) > 0$  contradicts the fact that $t \mapsto B(\phi(t))$ must be nonincreasing on $[t_1,t_2]$.  When $B$ is continuous,   we conclude that the (zero-sublevel) set $K$ is closed.   Hence,  there exists $0 \leq t_1 < t_2$ such that $\phi(t_2) \in  U(\partial K) \backslash K$,  $\phi(t_1) \in \partial K$,  and $\phi((t_1,t_2]) \subset U(\partial K) \backslash K$.   Now,  having   $t \mapsto B(\phi(t))$ nonincreasing on $(t_1,t_2]$ allows us to conclude,  using continuity of both $B$ and $\phi$,  that $t \mapsto B(\phi(t))$ is nonincreasing on 
 $[t_1,t_2]$.  But at the same time $B(\phi(t_1)) \leq 0$ and $B(\phi(t_2)) > 0$, which yields to a contradiction.  }
\end{proof}

Next, we recall from \cite{Sanfelice:monotonicity},  \cite{9705088},  \cite{della2021piecewise},  and \cite{bacciotti1999stability} characterizations of \ref{item:star1} using infinitesimal inequalities.

\begin{lemma} [Characterization of monotonicity] \label{lemA9}
Consider system $\Sigma$ with $F$ satisfying Assumption \ref{ass1} and a scalar function $B : \mathbb{R}^n \rightarrow \mathbb{R}$.   Given an open set $O \subset \mathbb{R}^n$,  we consider the monotonicity property:
\begin{enumerate} [label={($\star'$)},leftmargin=*]
\item \label{item:star2} 
Along every solution $\phi$  to $\Sigma$ satisfying $\phi(\dom \phi) \subset O$,  the map $t \mapsto B(\phi(t))$ is nonincreasing.
\hfill $\bullet$
\end{enumerate}

\begin{enumerate}
\item When $B$ is lower semicontinuous and $F$ is locally Lipschitz, \ref{item:star2} is satisfied if and only if 
\\
$ \langle \partial_P B(x) ,  \eta  \rangle \subset \mathbb{R}_{\leq 0} \qquad \forall \eta \in F(x), \quad \forall x \in O.  $

\item When $B$ is locally Lipschitz,  \ref{item:star2} is satisfied if
\\
$
\langle \partial_C B(x), \eta \rangle \subset  \mathbb{R}_{\leq 0}  \quad \forall \eta \in F(x), 
\quad   \forall x \in  O.
$

\item When $B$ is locally Lipschitz and regular,  \ref{item:star2} is satisfied if ~~ $ \langle \partial_C B(x), \eta \rangle \subset  \mathbb{R}_{\leq 0}  ~ \forall \eta \in F_B(x), 
~   \forall x \in  O,
$
\\
where $ F_B :  \mathbb{R}^n  \rightrightarrows \mathbb{R}^n$ is given by
\begin{align*}
\hspace{-0.4cm} F_B(x)  :=  \{ \eta \in F(x) : \text{Card} \left(\langle \partial_C B(x), \eta \rangle \right) = 1    \}, 
\end{align*}
where $\text{Card}\left(\langle \partial_C B(x), \eta \rangle \right)$ denotes the number of elements in the set $\left(\langle \partial_C B(x), \eta \rangle \right)$.

\item  When $B$ is locally Lipschitz and $F$ is continuous,  \ref{item:star2} is satisfied if and only if
\\
$ \langle \nabla B(x) , \eta  \rangle \leq 0 \quad  \forall \eta \in F(x), \quad \forall x \in O\backslash \Omega,  $ 
\\
where $\Omega \subset \mathbb{R}^n$ is the set of points where $B$ is not differentiable.  

\item When $B$ is continuously differentiable, \ref{item:star2} is satisfied if
$\langle \nabla B(x) , \eta  \rangle \leq 0 ~ \forall \eta \in F(x),  ~ \forall x \in O$. 

\item When $B$ is continuously differentiable and $F$ is continuous,  \ref{item:star2} is equivalent to the condition in the previous item.
\end{enumerate}
\end{lemma}

\begin{remark}
\blue{The first item is originally proved in \cite[Theorem 6.3]{clarke2008nonsmooth}.   
A proof of the second item can be found in 
\cite[Theorem 3]{Sanfelice:monotonicity}. 
The third item is originally proved in \cite{bacciotti1999stability}.    The proof of the fourth item can be found in \cite{della2021piecewise}.  The fifth item is a particular case of the second item.  
 Finally,    the sixth item is a particular case of the fourth.  
 }
\end{remark}

\ifitsdraft
We now recall Michael's Selection Theorem. 
\begin{lemma}[Selection Theorem \cite{michael1956continuous}] \label{lemA6}
Consider a closed set $K \subset \mathbb{R}^n$ and a set-valued map $\Phi : K \rightrightarrows \mathbb{R}^m$ that is lower semicontinuous with $\Phi(x)$ nonempty, closed, and convex for all $x \in K$. Then, for each closed set $A \subset K$ and for each $\phi : A \rightarrow \mathbb{R}^m$ a selection of $\Phi$ on $A$; namely, $\phi$ is continuous and 
$\phi(x) \in \Phi(x)$ for all $x \in A$,
there exists $\phi_e : K \rightarrow \mathbb{R}^m$ extending 
$\phi$ to $K$ that is a selection of $\Phi$ on $K$; namely, 
$\phi_e$ is continuous,  
$$ \phi_e(x) = \phi(x) \quad \forall x \in A \quad \text{and} \quad \phi_e(x) \in \Phi(x) \quad \forall x \in K. $$ 
\end{lemma}
\fi

\ifitsdraft

\red{
The following  result is in \cite[Theorem 4.3.4]{Aubin:1991:VT:120830}. 
\begin{lemma} [Contractivity] \label{lemA001}
Consider system $\Sigma$ such that Assumption \ref{ass1} holds. Let $K \subset \mathbb{R}^n$ be a closed set with nonempty interior and let $x_o \in \partial K$. If $F(x_o) \subset D_{\text{int}(K)}(x_o)$,
where $D_K$ is the \textit{Dubovitsky-Miliutin} cone of $K$ defined as 
\begin{align} \label{eq.cone2}
D_K(x) & :=  \{v \in \mathbb{R}^n : \exists \varepsilon >0: \nonumber 
\\ &
x + \delta (v + w) \in K~\forall \delta \in (0,\varepsilon],~\forall w \in \varepsilon \mathbb{B}  \}.
\end{align}
Then, for each solution $\phi$ to $\Sigma$ starting from $x_o$, there exists $T > 0$ such that $\phi((0,T],x_o) \subset \mbox{int}(K)$. 
\end{lemma} 

The following Lemma characterizes the cone $D_{\mbox{int}(K)}$ at the boundary of the set $K$ in terms of the barrier function candidate defining the set; see \cite[Lemma 4]{draftautomatica}.

\begin{lemma} [Characterization of $D_K$]  \label{lemA002}
Let a continuously differentiable function  $B : \mathbb{R}^n \rightarrow \mathbb{R}$ generating the set $K \subset \mathbb{R}^n$ according to \eqref{eq.4}.  Let $x \in \partial K$ such that $\nabla B(x) \neq 0$.  Then,   
\begin{center}
$D_{\mbox{int}(K)}(x) = \left\{ v \in \mathbb{R}^n : \langle \nabla B(x), v \rangle < 0 \right\}$. 
\end{center}
\end{lemma} }
 \fi

\begin{proposition}  \label{propagc}
\blue{Consider the sets $A \subset \mathbb{R}^m$ and $K \subset \mathbb{R}^n$, and let the system $\Sigma$ in \eqref{eqSigmaagc}. 

\begin{itemize}
\item If $\Sigma$ is safe with respect to $(K,\mathbb{R}^n \backslash K)$, then  $\Sigma_w$ weakly satisfies the contract $\mathcal{C} := (A,K)$;

\item If $K$ is closed, $f$ is locally Lipschitz, and $\Sigma$ is robustly safe with respect to $(K,\mathbb{R}^n \backslash K)$, then  $\Sigma_w$ strongly satisfies the contract $\mathcal{C} := (A,K)$,
\end{itemize}
where the system $\Sigma_w$ is defined in (\ref{eqn:sys_pert}).
}
\end{proposition}
\begin{proof}
\blue{
To prove the first item, we consider a maximal solution pair $(w,\phi) : \dom \phi \rightarrow \mathbb{R}^m \times \mathbb{R}^n$ to $\Sigma_u$. The first item of the weak contract's satisfaction is satisfied by definition of safety with respect to $(K,\mathbb{R}^n \backslash K)$. Now, select $t \in \dom \phi$ and assume that $w([0,t]) \subset A$. Hence, the solution pair 
$(w,\phi) : [0,t] \rightarrow \mathbb{R}^m 
\times \mathbb{R}^n$ is a solution to $\Sigma$, which is safe with respect to the set $(K,\mathbb{R}^n \backslash K)$. Hence, $\phi([0,t]) \subset K$, which implies weak satisfaction of the contract $\mathcal{C}$.
 
 Now, we assume that the system $\Sigma$ is robustly safe with respect to the set $(K,\mathbb{R}^n \backslash K)$, and let us show that the system $\Sigma_w$ strongly satisfies the contract $\mathcal{C}$. Consider a maximal solution pair  $(w,\phi) : \dom \phi \rightarrow \mathbb{R}^m\times \mathbb{R}^n$ to $\Sigma_w$. Consider $t \in dom \phi$ and  assume that $w([0,t]) \subset A$. 
 First, from the weak satisfaction of the contract, one gets that $\phi([0,t]) \in K$. Moreover, using the fact  $\Sigma$ is robustly safe with respect to  $(K,\mathbb{R}^n \backslash K)$, we conclude that $\phi(t) \in K \backslash \partial K = \mbox{int}(K)$, since, otherwise, using \cite{puri1995varepsilon}, we conclude that, starting from $\phi(t)$,  the solutions to $\Sigma_\epsilon$, with $\epsilon$ being a robustness margin,  would reach the exterior of $K$. Hence, using the fact that the set $\mbox{int}(K)$ is open, we conclude the existence of $\eta >0$ such that $\phi(t) + \eta \mathbb{B} \subseteq K$, which implies from local boundedness of $F$ and the existence of 
 non-trivial solutions to $\Sigma$ starting from $\phi(t)$, under continuity of $F$, that there exists  $\delta >0$ such that $\phi([t,t+\delta]) \subset K$. Hence, one gets that $\phi([0,t+\delta]) \subset K$, and $\Sigma_w$ strongly satisfies the contract $\mathcal{C} :=(A,K)$.}
 \end{proof}
\begin{lemma} \label{lemA4}
Consider a closed set $K \subset \mathbb{R}^n$ and an upper semicontinuous function
$f : K \rightarrow \mathbb{R}_{<0}$. 
Then, there exists a continuous function $g : K \rightarrow \mathbb{R}_{<0}$ such that  
$0 > g(x) \geq f(x)$ for all $x \in K$.  
\end{lemma}
\begin{proof}
We propose to grid the set $K$ using a sequence of nonempty compact subsets 
$\{K_i\}^{N}_{i=1}$, where $N \in \{1,2,...,\infty\}$. That is, we assume that 
$\bigcup^{N}_{i = 1} K_i = K$. Furthermore, for each $i \in \{1,2,...,N\}$, there exists a finite set 
$\mathcal{N}_i \subset \{1,2,...,N \}$ 
such that 
$K_i \cap K_j = \emptyset$ for all 
$j \notin \mathcal{N}_i$, and 
\begin{align} \label{eqch}
\mbox{int}_K (K_i \cap K_j) = \emptyset \qquad \forall j \in \mathcal{N}_i, 
\end{align}
where $\mbox{int}_K(\cdot)$ denotes the relative interior of $(\cdot)$ with respect to the closed set $K$. In words,  \eqref{eqch} implies that the sets $K_j$, $j \in \mathcal{N}_i$, intersect with the set $K_i$ without overlapping.  Such a decomposition always exists according to the Whitney covering lemma \cite{10.2307/1989708}.
Now, for each $i \in \{1,2,...,N\}$, we let 
$\alpha_i := \sup \{ f(x) : x \in K_i \}$. 
Using \cite[Theorem B.2]{puterman2014markov}, we conclude that, 
for each $i \in \{1,2,...,N\}$, there exists $x^*_i \in K_i$ such that $\alpha_i := f(x^*_i) < 0$. Next, we let 
$l : K \rightarrow \mathbb{R}_{<0}$ given by 
$l(x) := \max \{\alpha_i : 
i \in \{1,2,...,N\} ~ \text{s.t.} ~ x \in K_i \}$,
and we show that  $l(x) \geq f(x)$ for all $x \in K$.  Indeed, for $i \in \{1,2,...,N \}$, we distinguish two complementary situations. If $x \in \mbox{int}_{K}(K_i)$, the latter inequality holds for free. 
 Otherwise, there must exist 
$\mathcal{N}^x_i \subset \mathcal{N}_i$ such that $x \in \bigcap_{j \in \mathcal{N}^x_i} K_j \cap K_i$ and  $x \notin K_k \quad \forall k \neq \{i, \mathcal{N}^x_i\}. $
Thus, $\alpha_j \geq f(x)$ for all $j \in \{i,  \mathcal{N}^x_i \}$ and at the same time $l(x) = \max \{ \alpha_j : 
j \in \{i,  \mathcal{N}^x_i \}  \} \geq f(x)$. Furthermore, it is easy to see that the function $l$ is lower semicontinuous. Hence, we are able to define the set-valued map $G_{fl} : K \rightrightarrows \mathbb{R}_{< 0}$ given by 
$G_{fl} (x) := [f(x),l(x)]$. Note that $G_{fl}$ is lower semicontinuous with nonempty, closed, and convex images for all $x \in K$. Hence, using the selection theorem in \cite{michael1956continuous}, we conclude the existence a continuous function $g : K \rightarrow \mathbb{R}$ such that $0 > l(x) \geq  g(x) \geq f(x)$.
\end{proof}

\begin{lemma} \label{lemA3}
Consider a closed 
{subset} $K \subset \mathbb{R}^n$ and a set-valued map $G : K \rightrightarrows \mathbb{R}_{<0}$ that is upper semicontinuous with $G(x)$ nonempty, compact, and convex for all 
$x \in K$. Then, there exists  $g : K \rightarrow \mathbb{R}_{<0}$  continuous such that  
$0 > g(x) \geq \sup \{\eta: \eta \in G(x)\}$ for all $x \in K$.  
\end{lemma}
\begin{proof}
Let the function 
$f : K \rightarrow \mathbb{R}_{<0}$ given by
$f(x) := \sup \{\eta: \eta \in G(x)\}$. 
Using \cite[Theorem 1.4.16]{aubin2009set} we conclude that $f$ is upper semicontinuous. Finally, using Lemma \ref{lemA4}, we conclude that there exists $g : K \rightarrow \mathbb{R}_{<0}$ continuous such that $g(x) \geq f(x)$ for all $x \in K$.
\end{proof}

\bibliography{biblio}
\bibliographystyle{ieeetr}     

\end{document}